\documentstyle[12pt]{article}

\begin{document}
\def\comment#1{\marginpar{{\scriptsize #1}}}
\def\framew#1{\fbox{#1}}
\def\framep#1{\noindent \fbox{\vbox{#1}}}
\def\framef#1{\fbox{\vbox{ #1 }}}
\def\be{\begin{equation}}
\def\bea{\begin{eqnarray}}
\def\ee{\end{equation}}
\def\eea{\end{eqnarray}}
\newcommand{\sect}[1]{\setcounter{equation}{0}\section{#1}}
\newcommand{\subsect}[1]{\subsection{#1}}
\newcommand{\subsubsect}[1]{\subsubsection{#1}}
\renewcommand{\theequation}{\thesection.\arabic{equation}}
\renewcommand{\thefootnote}{\fnsymbol{footnote}}
\def\fn{\footnote}
\footskip 1.0cm
\def\sxn#1{\bigskip\medskip \sect{#1} \smallskip}
\def\subsxn#1{\medskip \subsect{#1} \smallskip}
\def\subsubsxn#1{\medskip \subsubsect{#1} \smallskip}
\newtheorem{proposition}{Proposition}[section]  
\def\bprop{\medskip\begin{proposition}~~~\smallskip \it}
\def\eprop{\end{proposition}\bigskip}
\def\proof{\bigskip \noindent {\it Proof.} \ }
\newtheorem{naming}{Definition}[section]   
\def\bnam{\medskip\begin{naming}~~~\smallskip \it}
\def\enam{end{naming}\bigskip}
\newtheorem{example}{Example}[section]   
\def\bexam{\medskip\begin{example} ~~\\ \rm}
\def\eexam{\end{example}\bigskip}
\bibliographystyle{unsrt}
\def\br{\begin{thebibliography}{abc}}
\def\er{\end{thebibliography}}
\def\rf{\bibitem}

%
%
\def\pinco{| \! / \! |}
\def\cstars{$C^*$-algebras }
\def\cstar{$C^*$-algebra }
\def\unit{I\!\!I}
\def\norm#1{\parallel #1 \parallel}
\def\abs#1{\left| #1\right|}
\def\ha{\widehat{\cal A}}
\def\hc{\widehat{\cal C}}
\def\prim{$Prim \ca$ }

\def\ca{{\cal A}}
\def\cb{{\cal B}}
\def\cc{{\cal C}}
\def\cd{{\cal D}}
\def\ce{{\cal E}}
\def\cf{{\cal F}}
\def\cg{{\cal G}}
\def\ch{{\cal H}}
\def\ci{{\cal I}}
\def\cj{{\cal J}}
\def\ck{{\cal K}}
\def\cl{{\cal L}}
\def\cm{{\cal M}}
\def\cn{{\cal N}}
\def\co{{\cal O}}
\def\cp{{\cal P}}
\def\cq{{\cal Q}}
\def\czr{{\cal R}}
\def\cs{{\cal S}}
\def\ct{{\cal T}}
\def\cu{{\cal U}}
\def\cv{{\cal V}}
\def\cw{{\cal W}}
\def\cx{{\cal X}}
\def\cy{{\cal }}
\def\cz{{\cal Z}}
\def\bar#1{\overline{#1}}
\def\what{\widehat}
\def\wtilde{\wtilde}
\def\bra#1{\left\langle #1\right|}
\def\ket#1{\left| #1\right\rangle}
\def\EV#1#2{\left\langle #1\vert #2\right\rangle}
\def\VEV#1{\left\langle #1\right\rangle}
%
%
\def\sdp{\hbox{ \raisebox{.25ex}{\tiny $|$}\hspace{.2ex}{$\!\!\times $}} }
\def\pa{\partial}
\def\del{\nabla}
\def\a{\alpha}
\def\b{\beta}
\def\c{\raisebox{.4ex}{$\chi$}}
\def\d{\delta}
\def\e{\epsilon}
\def\f{\phi}
\def\g{\gamma}
\def\h{\eta}
\def\i{\iota}
\def\j{\psi}
\def\k{\kappa}
\def\l{\lambda}
\def\m{\mu}
\def\n{\nu}
\def\o{\omega}
\def\p{\pi}
\def\q{\theta}
\def\r{\rho}
\def\s{\sigma}
\def\t{\tau}
\def\u{\upsilon}
\def\x{\xi}
\def\z{\zeta}
\def\D{\Delta}
\def\F{\Phi}
\def\G{\Gamma}
\def\J{\Psi}
\def\L{\Lambda}
\def\O{\Omega}
\def\P{\Pi}
\def\Q{\Theta}
\def\S{\Sigma}
\def\U{\Upsilon}
\def\X{\Xi}
\def\Z{\Zeta}
%

%
%
%
\def\inbar{\,\vrule height1.5ex width.4pt depth0pt}
\def\IC{\relax\,\hbox{$\inbar\kern-.3em{\rm C}$}}
\def\ID{\relax{\rm I\kern-.18em D}}
\def\IF{\relax{\rm I\kern-.18em F}}
\def\IH{\relax{\rm I\kern-.18em H}}
\def\II{\relax{\rm I\kern-.17em I}}
\def\I1{\relax{\rm 1\kern-.28em l}}
\def\IM{\relax{\rm I\kern-.18em M}}
\def\IN{\relax{\rm I\kern-.18em N}}
\def\IP{\relax{\rm I\kern-.18em P}}
\def\IQ{\relax\,\hbox{$\inbar\kern-.3em{\rm Q}$}}
\def\IZ{\relax\,\hbox{$\inbar\kern-.3em{\rm Z}$}}
\def\IR{\relax{\rm I\kern-.18em R}}
\font\cmss=cmss10 \font\cmsss=cmss10 at 7pt
\def\Z{\relax\ifmmode\mathchoice
{\hbox{\cmss Z\kern-.4em Z}}{\hbox{\cmss Z\kern-.4em Z}}
{\lower.9pt\hbox{\cmsss Z\kern-.4em Z}}
{\lower1.2pt\hbox{\cmsss Z\kern-.4em Z}}\else{\cmss Z\kern-.4emZ}\fi}
%
%
\def\bc{{\bf C}}
\def\br{{\bf R}}
\def\bz{{\bf Z}}
\def\bn{{\bf N}}
\def\bm{{\bf M}}
\def\Up{\Uparrow}
\def\up{\uparrow}
\def\Dn{\Downarrow}
\def\dn{\downarrow}
\def\Ra{\Rightarrow}
\def\ra{\rightarrow}
\def\La{\Leftarrow}
\def\la{\leftarrow}
\def\iff{\Leftrightarrow}
\def\Up{\Uparrow}
\def\up{\uparrow}
\def\Dn{\Downarrow}
\def\dn{\downarrow}
\def\Rt{\Rightarrow}
\def\rt{\rightarrow}
\def\Lt{\Leftarrow}
\def\lt{\leftarrow}
\thispagestyle{empty}
\setcounter{page}{0}

\hfill \today

\vskip2cm

\centerline {\Large NONCOMMUTATIVE LATTICES}
\vspace{3mm}
\centerline {\Large AND THE ALGEBRAS OF}
\vspace{3mm}
\centerline{\Large THEIR CONTINUOUS FUNCTIONS}
\vspace{0.5cm}
\centerline {Elisa Ercolessi$^1$,
             Giovanni Landi$^{2,3}$,
             Paulo Teotonio-Sobrinho$^{2,4}$}
\vspace{0.5cm}
\centerline {\it $^1$ Dipartimento di Fisica, Universit\`a di
Bologna and INFM,}
\centerline{\it Via Irnerio 46, I-40126, Bologna, Italy.}
\vspace{1mm}
\centerline {\it $^2$ The E. Schr\"odinger International Institute for
Mathematical Physics,}
\centerline{\it Pasteurgasse 6/7, A-1090 Wien, Austria.}
\vspace{1mm}
\centerline{\it $^3$ Dipartimento di Scienze Matematiche,
Universit\`a di Trieste,}
\centerline{\it P.le Europa 1, I-34127, Trieste, Italy}
\centerline {\it and INFN, Sezione di Napoli, Napoli, Italy.}
\vspace{1mm}
\centerline {\it $^4$ Department of Physics, University of Illinois at
Chicago,} 
\centerline {\it 60607-7059 Chicago, IL, USA}
\centerline {\it and Universidade de Sao Paulo, Instituto de Fisica - DFMA,}
\centerline{\it Caixa Postal 66318, 05389-970, Sao Paulo, SP,
Brasil\footnote{Permanent address}.} 
\vspace{.3cm}

\begin{abstract}
Recently a new kind of approximation to continuum topological
spaces has been introduced, the approximating spaces being partially 
ordered sets (posets) with a finite or at most a countable number of 
points. The partial order endows a poset with a nontrivial 
non-Hausdorff topology. Their ability to reproduce important 
topological information of the continuum has been the main 
motivation for their use in quantum physics. Posets are truly 
noncommutative spaces, or {\it noncommutative lattices}, since they 
can be realized as structure spaces of noncommutative $C^*$-algebras. 
These noncommutative algebras play the same role of the algebra of 
continuous functions ${\cal C}(M)$ on a Hausdorff topological space 
$M$ and can be thought of as algebras of operator valued functions on 
posets. In this article, we will review some mathematical results that 
establish a duality between finite posets and a certain class of 
C$^*$-algebras. We will see that the algebras in question are all 
postliminal approximately finite dimensional (AF) algebras.

\end{abstract}

\newpage

\setcounter{footnote}{0}

\sxn{Introduction}\label{se:int}

It is well known that the standard discretization methods used in quantum
physics (where a manifold is replaced by a lattice of points with the discrete
topology) are not able to describe any significant topological attribute of the
continuum, this being equally the case for both the local and global
properties. 
For example, there is no nontrivial concept of winding number 
and hence no way to formulate theories with topological solitons or
instantons on 
these lattices.

A new kind of finite approximation to continuum topological
spaces has been first introduced in \cite{So}, with the name of posets or
partially ordered sets. As we will see in section \ref{se:nc}, posets are also
$T_0$ topological spaces and can reproduce
important topological properties of the continuum, such as the homology and the
homotopy groups, with remarkable fidelity \cite{So,Al}.  This ability
to capture 
topological information has been the main motivation for their use in quantum
physics in substitution of the ordinary discrete lattices. In
\cite{BBET}, quantum 
mechanics has been formulated on posets and it has been proved that it
is possible 
to study nontrivial topological configurations, such as $\theta$-states for
particles on the poset approximations to a circle. Some promising
results have also 
been obtained in the formulation of solitonic field theories
\cite{BBET} as well as 
of gauge field theories \cite{Field}.

In \cite{lisbon}, the poset approximation scheme has been developed
in a novel direction. Indeed, it has been observed that posets are truly
noncommutative spaces, or {\it noncommutative lattices}, since they
can be realized 
as structure spaces (spaces of irreducible representations) of noncommutative
$C^*$-algebras. These noncommutative algebras play the same r\^{o}le
of the algebra 
of continuous functions $\cc(M)$ on a manifold $M$ and can be thought of as
algebras of operator valued functions on posets. This naturally leads
to the use of 
noncommutative geometry \cite{Co} (see also \cite{VG}) as the tool to rewrite
quantum theories on posets and gives a remarkable connection between
topologically 
meaningful finite approximations to quantum physics and noncommutative
geometry. 

The duality relation between Hausdorff topological spaces and
commutative C$^*$-algebras is provided by the Gel'fand-Naimark
theorem. There is no analogue of this theorem in the noncommutative
setting. In this article, we will review how it is possible
to establish a relation between finite posets and a particular class of
noncommutative C$^*$-algebras. For such class of algebras the situation is very
similar to the commutative case. We will see that the algebras in
question are all approximately finite dimensional (AF) postliminal algebras 
\cite{Br,Go}, i.e. \cstars that can be approximated in norm by a
sequence of finite 
dimensional matrix algebras and whose irreducible representations are 
completely characterized by the kernels. This is exactly what makes
them of some 
interest in mathematics: they present virtually all the attributes and
complications of other infinite dimensional algebras, but many techniques and
results valid in the finite dimensional case can be used in their
study. Thus, for 
example, a complete classification of AF \cstars is available \cite{Ef}.

These algebras have been first  extensively studied by
Bratteli \cite{Br}, who also introduced a diagrammatic representation 
which is very useful for the study of their algebraic properties. In particular
we will see how to use the Bratteli diagram of an AF algebra to construct its
structure space. Then we will see how, given any finite poset $P$, it
is possible 
to construct the Bratteli diagram of an AF algebra whose structure space is the
given poset $P$. Being noncommutative, this AF \cstars is far from
being unique. 
Indeed there is a whole family of AF algebras that have $P$ as structure space
and that can be classified by means of results due to Behncke 
and Leptin \cite{BL}.

In this article, we will not present  the classification of AF
$C^*$-algebras that can be formulated in terms of their algebraic K-theory
\cite{Ef}.  In view of their relation with posets, this would represent also a
first step in the construction of bundles and characteristic classes over
noncommutative lattices. This is indeed the content of
\cite{kt}, to which we refer the interested reader for a detailed analysis of
the K-theory of AF algebras.

This article is organized as follows. In section \ref{se:css} we review some
elementary algebraic concepts as well as the Gel'fand Naimark theorem
in order to 
clarify the notation and the terminology used in the sequel. In section
\ref{se:nc} we  briefly describe the topological approximation of continuous
spaces that leads to partially ordered sets (posets). In section
\ref{se:af}, AF 
\cstars are introduced and the connection between Bratteli diagrams and posets
is discussed in detail. Finally, in section \ref{se:bl} we present the 
classification theorem of the AF \cstars that have a poset as structure space
due to Behncke and Leptin. Several interesting examples will be examined
throughout the article.

\sxn{\cstars and Structure Spaces}\label{se:css}

Let us start with some elementary algebraic preliminaries \cite{Di,FD}
that will 
be also useful to set the notation. 

In the sequel, $\ca$ will always denote a \cstar over the  field of complex
numbers $\IC$. We remind that this means that $\ca$ is equipped with a norm of
algebra $|| \cdot || :\ca \rightarrow \IC$ (with respect to which $\ca$ is
complete) and an involution $^*: \ca \rightarrow \ca$, satisfying the identity:
\be
\norm{a^*a} = {\norm a}^2~,~~~ \forall~ a \in \ca~. \label{ss2}
\ee

The following are examples of commutative and noncommutative \cstars which
will be used in the article:
\begin{itemize}
\item[1)]
the (noncommutative) algebra $\IM(n,\IC)$ of $n\times n$ complex matrices $T$,
with $T^*$ given by the hermitian conjugate of $T$ and the squared norm
$\parallel T\parallel^2$ being equal to the largest eigenvalue of $T^* T$;
\item[2)]
the (noncommutative) algebra $\cb(\ch)$ of bounded operators $B$ on a separable
(infinite-dimensional)  Hilbert space $\ch$ as well as its subalgebra
$\ck(\ch)$ of compact operators. Now $^*$ is given by the adjoint and
the norm is 
the operator norm: $||B || = \sup_{||\xi || \leq 1} ||B \xi|| \; (\xi \in
\ch)$.

 \item[3)]
the (commutative) algebra $\cc(M)$ of continuous functions on a compact
Hausdorff topological space $M$, with * denoting complex conjugation and the
norm given by the supremum norm, $\norm f _\infty = \sup_{x\in M}|f(x)|$.
If $M$ is not compact but only locally compact, then one can consider
the algebra $\cc_0(M)$ of functions vanishing at infinity.
\end{itemize}

Notice that $\ck(\ch)$ and $\cc_0(M)$ (with $M$ only locally compact) are
examples of \cstars without unit $\unit$, contrary to $\IM(n,\IC)$ and
$\cb(\ch)$.

\subsxn{Commutative \cstars: The Gel'fand-Naimark
\mbox{Theorem}}\label{sub:gnt} 

In the third example above we have seen how it
is possible to associate a commutative \cstar with (without) unit, namely
$\cc(M)$ ($\cc_0(M)$), to any Hausdorff compact (locally compact) topological
space $M$. Viceversa, given {\it any} commutative \cstar $\cc$ with (without)
unit, it is possible to construct a Hausdorff compact (locally compact)
topological space $M$ such that $\cc$ is isometrically $*$-isomorphic to the
the algebra of continuous functions $\cc(M)$ ($\cc_0(M)$). This is precisely
the content of the Gel'fand-Naimark theorem \cite{FD} that will be
discussed in this 
paragraph. For simplicity we will consider only the case when $\cc$ is a
commutative \cstar with unit.

Given such a $\cc$, we let $\hc$ denote the {\it structure space} of
$\cc$, namely the space of equivalence classes of irreducible
$^*$-representations 
(IRR's) of $\cc$. The trivial IRR  given by $\cc\rightarrow \{0\}$ is
not included in $M$ and will therefore be ignored here and hereafter.
Since the \cstar $\cc$ is commutative, every IRR is one-dimensional, i.e. is
a (non-zero) linear  functional $\f : \cc \rt \IC$ satisfying $\f(ab) =
\f(a)\f(b)$ and $\f(a^*) = \overline{\f(a)}$, $\forall a, b \in
\cc$. It follows 
that $\f(\unit) = 1, \forall \f \in \hc$.
The space $\hc$ is made into a topological space by endowing it with the {\it
Gel'fand topology}, namely with the topology of pointwise convergence on
$\cc$. Then $\cc$ can be proved to be a Hausdorff compact topological space.

For a commutative $C^*$-algebra, two-irreducible representations are unitarily
equivalent if and only if they have the same kernel. Thus one can consider
also the space of kernels of IRR's, the so called {\it primitive spectrum}
$Prim \cc$. Now, these kernels are maximal ideals of
$\cc$  and, viceversa, any maximal ideal is the kernel of an
irreducible representation \cite{FD}.
Indeed, suppose that $\f \in \hc$, then $Ker(\f)$ is of codimension
$1$ and so is a maximal ideal of $\cc$. Conversely, suppose that $I$
is a maximal ideal of $\cc$, then the natural representation of $\cc$
on $\cc / I$ is irreducible, hence one-dimensional. It follows that
$\cc / I \cong \IC$, so that the quotient homomorphism
$\cc \rt \cc / I$ can be identified with an element $\f \in \hc$.
Clearly, $I = Ker(\f)$. Thus $Prim \cc$ can be thought of as the space of
maximal ideals.  As such, $Prim \cc$ is equipped with the {\it Jacobson} or
{\it hull kernel topology}, that will be described in the next
paragraph.

The map that to each class of unitary representations associates its kernel
gives a map $\;\hc \rt Prim \cc$, which turns out to be a homeomorphism of the
two topological spaces so that we can equivalently talk of the structure
space or of the primitive spectrum of a commutative $C^*$-algebra.

If $c \in \cc$, the {\it Gel'fand transform} $\hat{c}$ of $c$ is the
complex-valued function on $\hc$ given by
\be\label{gn1}
\hat{c}(\f) = \f(c)\; , \; \; \forall \f \in \hc\; .
\ee
It is clear that $\hat{c}$ is continuous for each $c$. We thus get the
interpretation of elements in $\cc$ as {\bf C}-valued continuous functions on
$\hc$. The Gel'fand-Naimark theorem states that all continuous functions on
$\hc$ are indeed of the form (\ref{gn1}) for some $c \in \cc$ \cite{FD}:
\bprop\label{th:gn}
Let $\cc$ be a commutative $C^*$-algebra. Then the
Gel'fand transform map $c  \mapsto \hat{c}$ is an isometric $*$-isomorphism of
the \cstar $\cc$ onto  the \cstar $\cc(\hc)$ (equipped with the supremum norm
$\norm{\cdot}_{\infty}$).
\eprop

Suppose now that $M$ is a compact topological space.
We have a natural $C^*$-algebra, $\cc(M)$, associated to it. It is
then natural to 
ask what is the relation between the Gel'fand space
$\what{\cc(M)}$ and $M$ itself. It turns out that this two spaces can be
identified both setwise and topologically.
We notice first that each $m \in M$ gives a complex homomorphism $\f_m \in
\what{\cc(M)}$ through the evaluation map:
\be\label{gn3}
\f_m : \cc(M) \rt \IC~, ~~~ \f_m(f) = f(m)~.
\ee
Let $I_m$ denote the kernel of $\f_m$, namely the maximal ideal of
$\cc(M)$ consisting of all functions vanishing at $m$. We have the
following theorem \cite{FD}:
\bprop\label{th:gn1}
The map $\Phi:m \mapsto \f_m$ given by (\ref{gn3}) is a homeomorphism of
$M$ onto  $\what{\cc(M)}$, namely $M \cong \what{\cc(M)}$. Moreover,
every  maximal ideal of $\cc(M)$ is of the form $I_m$ for some $m \in M$.
\eprop

In conclusion, the previous two theorems set up a one-to-one correspondence
between  the $*$-isomorphism classes of commutative \cstars and the
homeomorphism classes of locally compact Hausdorff spaces. If $\cc$ has a unit,
then $\hat{\cc}$ and $Prim\cc$ are compact.

\subsxn{Noncommutative Algebras and Associated Spaces}\label{se:nca}

The scheme described in the previous section cannot be directly generalized
to noncommutative $C^*$-algebras. There is more than one candidate for
the analogue of the topological space $M$. In particular, since
non-equivalent unitary transformations may now have the same kernel, we
have to distinguish even setwise between:
\begin{itemize}
\item[1)]
the {\it structure space} $\ha$ of $\ca$ or the space of all
unitary equivalence classes of irreducible $^*$-representations;
\item[2)]
the {\it primitive spectrum} \prim of $\ca$ or the space of
kernels of irreducible $^*$-representations. Any element of  \prim is
automatically a two-sided $^*$-ideal of $\ca$.
\end{itemize}

One can define a natural topology on both $\ha$ and $Prim\ca$.
While for a commutative \cstar the resulting topological spaces
are homeomorphic, this is no longer true in the noncommutative case. For
instance, in section (\ref{sub:dia}) we will describe a \cstar $\ca$ associated
to the Penrose tiling of the plane \cite{Co}, whose structure space
$\ha$ consists 
of an infinite countable set of points, whereas \prim consists of a
single point. 

The topology one puts on $\ha$ is called {\it regional topology}
\cite{FD} and is 
a generalization of the pointwise convergence we have used in the previous
paragraph, to which it reduces in the commutative case. This topology is
constructed by defining a basis of neighborhoods for the points
(classes of representations) of $\ha$ as follows. Given a $T\in \ha$, let us
denote with $\ch_T$ the Hilbert space of the representation $T$. Then an
open neighborhood of $T$ is identified by a finite sequence
$\xi_1,\xi_2,\cdots,\xi_n$ of vectors in $\ch_T$, a positive number $\e$
and a finite not void set $\cf\subset \ca$ by means of:
\bea
& & U(T;\e;\xi_1,\xi_2,\cdots,\xi_n;\cf) =: \{ T'\in \ha \, : \,
\exists \, \xi_1',\xi_2',\cdots,\xi_n'\in \ch_{T'} \mbox{ with } \nonumber \\
& & |(\xi'_i,\xi'_j)_{\ch_{T'}} - (\xi_i,\xi_j)_{\ch_{T}} | < \e \; , \;
|(T'(a)\xi'_i,\xi'_j)_{\ch_{T'}} - (T(a)\xi_i,\xi_j)_{\ch_{T}} | < \e
\nonumber \\ & & \mbox{ for } i,j=1,2,\cdots,n
\mbox{ and } \forall a\in \cf \}\; .
\eea

On \prim we instead define a closure operation as follows \cite{Di,FD}.
Given any subset $W$ of $Prim\ca$, the closure $\bar{W}$ of $W$ is by
definition 
the set of all elements in \prim containing the intersection of the elements
of $W$, namely
\be
\bar{W} = \{ I \in Prim{\ca} : \bigcap W \subset I\}~. \label{cl}
\ee
This `closure operation' satisfies the Kuratowski axioms \cite{Ke} and thus
defines a topology on $Prim\ca$, which is called {\it Jacobson} or {\it
hull-kernel topology}. With respect to this topology we have:
\bprop\label{pr1}
Let $W$ be a subset of $Prim\ca$. Then $W$ is closed if and only if
$W$ is exactly the set of primitive ideals containing some subset of $\ca$.

\proof
If $W$ is closed, by \ref{cl}, $W$ is the set of primitive ideals containing
$\bigcap W \subset I$. Conversely, let $V \subseteq \ca$. If $W$ is the set
of primitive ideals of $\ca$ containing $V$, then
$V \subseteq \bigcap W \subset I$, for all $I \in W$, so that $\bar{W}
\subset W$, 
and $\bar{W}= W$.
\eprop

In general $\ha$ and \prim fail to be Hausdorff (or $T_2$).
Recall \cite{Ke} that a topological space is called $T_0$ if for any
two distinct 
points of the space there is a neighborhood of one of the two points
which does not 
contain the other. It is called $T_1$ if any point of the space is closed.
It is called $T_2$ if there exist disjoint neighborhoods of any two
points. Whereas 
nothing can be said on the separability properties of $\ha$, it turns out that
\prim is always a $T_0$ space and that it is $T_1$ if and only if  all
primitive 
ideals in $\ca$ are maximal, as it is established by the following propositions
\cite{Di,FD}.

\bprop\label{pr2}
The space \prim is a $T_0$ space.

\proof
Suppose $I_1$ and $I_2$ are two distinct points of \prim so that say $I_1
\not\subset I_2$. Then the set of those $I \in$ \prim which contain $I_1$ is
a closed subset $W$ (by proposition \ref{pr1}) such that $I_1 \in W$ and $I_2
\not\in W$. Then its complement $W^c$ is an open set containing 
$I_2$ and not $I_1$. 
\eprop

\bprop\label{pr3}
Let $I \in $ $Prim\ca$. Then the point $\{I\}$ is closed in \prim if and
only if $I$ is maximal among primitive ideals.

\proof
Indeed $\bar{\{I\}}$ is just the set of primitive ideals of
$\ca$ containing $I$.
\eprop

One should be aware of the fact that, for a general algebra $\ca$, a
maximal ideal 
needs not be primitive. This is however always the case if $\ca$ admits a unit
\cite{Di}.

As in the commutative case, both $\ha$ and \prim are locally compact
topological spaces. In addition, if $\ca$ has a unit, then they are compact.
Notice that, in general, $\ha$ being compact does not implies that $\ca$ has
a unit. For instance, the algebra $\ck(\ch)$ of compact operators
on an infinite dimensional Hilbert space $\ch$ has no unit but its
structure space consists of a single point.

Let us now come to a comparison between the space $\ha$ and
$Prim\ca$. There is a 
canonical surjection of $\ha$ onto $Prim\ca$, given by the map that to each
IRR $\pi$ associates its kernel $ker \pi$. The pull-back of the Jacobson
topology from \prim to $\ha$ defines another topology on the latter that
turns out to be equivalent to the regional topology defined above
\cite{FD}. But 
$\ha$ and \prim are homeomorphic only under the hypotheses stated below
\cite{FD}.

\bprop
Let $\ca$ be a \cstar, then the following conditions are equivalent:
\begin{itemize}
\item[(i)] $\ha$ is a $T_0$ space.
\item[(ii)] Two irreducible representations of $\ha$ with the same
kernel are equivalent.
\item[(iii)] The canonical map $\ha \longrightarrow$ \prim is a homeomorphism.
\end{itemize}

\proof
By construction, a subset $\cs \in \ha$ will be closed if and only if it is of
the form $\{\pi \in \ha : ker \pi \in W \}$ for some $W$ closed in \prim. As a
consequence, given any two (classes of) representations $\pi_1, \pi_2
\in \ha$,  the representation $\pi_1$ will be in the closure of $\pi_2$ if and
only if $ker \pi_1$ is in the closure of $ker \pi_2$, or, by
proposition \ref{pr3} 
if and only if $ker \pi_2 \subset ker \pi_1$. In turn, $\pi_1$ and $\pi_2$ are
one in the closure of the other if and only if $ker \pi_2 = ker \pi_1$.
Therefore, $\pi_1$ and $\pi_2$ will not be distinguished by the topology of
$\ha$ if and only if they have the same kernel. On the other side, if
$\ha$ is $T_0$ one is able to distinguish points. It follows that $(i)$ implies
$(ii)$, namely that if $\hat{\ca}$ is a $T_0$ space, two
representations with the 
same kernel must be equivalent.
The other implications are obvious. 
\eprop

\sxn{Noncommutative Lattices}\label{se:nc}

For convenience, we will review 
in this section the content of \cite{So,Al,BBET}, where it is shown how it is
possible to approximate a continuum topological space by means of a finite or
countable set of points $P$ \cite{So} which, being equipped with a
partial order 
relation, is a \underline{p}artially \underline{o}rdered \underline{set}
or a {\it poset }. As explained there, these approximating 
spaces are able to reproduce important topological properties of the 
continuum. Moreover, in section \ref{sub:dia} we will see that any of 
these spaces can be identified with the space $\ha = Prim\ca$ of 
primitive ideals of some (noncommutative) AF algebra $\ca$, which thus 
plays the r\^{o}le of the algebra
of continuous functions on $P$ \cite{lisbon}. This fact will make any poset a
truly noncommutative space \cite{Co}, hence also the name {\it noncommutative
lattice}.

This is the reason why, in this article, we will consider only a
special class of 
algebras, namely postliminal approximately finite (AF) algebras. In section
\ref{se:af} we will see in detail that AF algebras are approximated in norm by
direct sums of finite dimensional matrices \cite{Br,Go}. As for postliminal we
refer to \cite{Di} for the exact definition. For what concerns this article, we
need just to know that, as a consequence of general theorems, this implies that
$\ha$ and \prim are homeomorphic. In other words, in the following we
will have to 
deal only with structure spaces (or primitive spectrum spaces) which are $T_0$
locally compact topological spaces.

\subsection{The Finite Topological Approximation}\label{sub:fta}

Let $M$ be a continuum topological space. Experiments are never so
accurate that they can detect events associated with points of $M$, rather they
only detect events as occurring in certain sets $O_\l$. It is therefore natural
to identify any two points $x$, $y$ of $M$ if they can never be separated or
distinguished by the sets $O_\l$.

Let us assume that each $O_\l$ is open and that the family $\{O_\l\}$
covers $M$:
\be
   M=\bigcup _\l O_\l~. \label{2.1}
\ee
We also assume that $\{O_\l \}$ is a topology for $M$ \cite{Ke}.
This implies that both $O_\l\cup O_\m$ and
$O_\l\cap O_\m$ are in $\cu$ if $O_{\l,\m}\in \cu$.
This hypothesis is physically consistent because experiments can
isolate events in $O_\l\cup O_\m$ and $O_\l\cap O_\m$ if they can do so
in $O_\l$ and $O_\m$ separately, the former by detecting an event in either
$O_\l$ or $O_\m$, and the latter by detecting it in both $O_\l$ and $O_\m$.

Given $x$ and $y$ in $M$, we write $x\sim y$ if every set $O_\l$
containing either point $x$ or $y$ contains the other too:
\be
   x\sim y \mbox{  means  } x\in O_\l \iff y\in O_\l
~~~\mbox{for every}~~ O_\l~ . \label{2.2}
\ee
Then $\sim $ is an equivalence relation, and it is reasonable to replace $M$
by $M\, / \sim \equiv P(M)$ to reflect the coarseness of observations.

We assume that the number of sets $O_\l$ is finite when $M$ is compact so that
$P(M)$ is an approximation to $M$ by a finite set in this case. When $M$ is
not compact, we assume instead that each point has a neighborhood
intersected by 
only finitely many $O_\l$,  so that $P(M)$ is a ``finitary"
approximation to $M$ 
\cite{So}. In the notation we employ, if $P(M)$ has $N$ points, we sometimes
denote it by $P_N(M)$.

The space $P(M)$ inherits the quotient topology from $M$ \cite{Ke}, i.e.
a set in $P(M)$ is declared to be open if its inverse image for
$\Phi$ is open in $M$, $\Phi$ being the map from $M$ to $P(M)$ obtained by
identifying equivalent points. The topology generated by these open sets is
the finest one compatible with the continuity of $\Phi$.

Let us illustrate these considerations for a cover of $M=S^1$ by four open
sets as in Figure \ref{fi:circle}(a).
In that figure, $O_{1},O_{3}\subset O_2\cap O_4$.
Figure \ref{fi:circle}(b) shows the corresponding discrete space
$P_4(S^1)$, the 
points $x_i$ being images of sets in $S^1$. The map $\Phi\, :\, S^1\rightarrow
P_4(S^1)$ is given by
$$
  O_1 \rightarrow x_1, ~ ~ ~ ~ ~
        O_2\setminus [O_2\cap O_4] \rightarrow x_2~,
$$
\be
  O_3 \rightarrow x_3, ~ ~ ~ ~ ~
        O_4\setminus [O_2\cap O_4] \rightarrow x_4 ~ .\label{2.3a}
\ee
The quotient topology for $P_4(S^1)$ can be read off from Figure
  \ref{fi:circle}, 
the open sets being
\be
\{x_1\}~,  ~ \{x_3\}~,  ~\{x_1,x_2,x_3\}~,  ~ \{x_1,x_4,x_3\}~,
 \label{2.4}
\ee
and their unions and intersections (an arbitrary number of the latter being
allowed as $P_4(S^1)$ is finite).

\begin{figure}[htb]
\begin{center}
\begin{picture}(200,100)(10,-40)
\put(10,10){\circle{40}}
\put(10,5){\oval(60,70)[t]}
\put(10,10){\oval(50,60)[b]}
\put(10,50){$O_2$}
\put(10,-30){$O_4$}
\put(-40, 7){$O_1$}
\put(50, 7){$O_3$}
\put(115,30){$\Phi$}
\put(100,10){\vector(1,0){40}}
\put(200,10){\circle*{4}}
\put(250,10){\circle*{4}}
\put(225,35){\circle*{4}}
\put(225,-15){\circle*{4}}
\put(185,10){$x_1$}
\put(255,10){$x_3$}
\put(225,40){$x_2$}
\put(225,-25){$x_4$}
\end{picture}
\end{center}
\caption{\label{fi:circle}}
\centerline{{\footnotesize  The covering of $S^1$ that gives rise to the
poset $P_4 (S^1 )$.}}
\vskip1cm
\end{figure}

Notice that our assumptions allow us to isolate events in certain sets of the
form $O_\l\setminus [O_\l\cap O_\m]$ which may not be open. This means that
there are in general points in $P(M)$ coming from sets which are not open in
$M$ and therefore are not open in the quotient topology. This implies that
in general $P(M)$ is neither Hausdorff nor $T_1$.
However, it can be shown \cite{So} that it is always a $T_0$ space. For
example, given the points $x_1$ and $x_2$ of $P_4(S^1)$, the open set $\{x_1
\}$ contains $x_1$ and not $x_2$, but there is no open set containing $x_2$
and not $x_1$.\\

We will see now how the topological properties of $P(M)$ can also be encoded
in a combinatorial structure, namely a partial order relation, that
can be defined 
on it.

Since $P(M)$ is finite (finitary), its topology is generated by the
smallest open neighborhoods $O_x$ of its points $x$. It is possible to
introduce a partial order relation $\preceq $ \cite{Al,St} by declaring
that:
\be
x\preceq y  \; \Leftrightarrow O_x \;  \subset O_y \; . \label{order}
\ee
In this way $P(M)$ becomes a \underline{p}artially \underline{o}rdered
\underline{set} or a {\it{poset}}.

Later, we will write $x\prec y$ to indicate that $x\preceq y$ and
$x\neq y$. A point 
$x \in P$ such that there exists no $y \in P$ with $x \prec y$ ($x \succ y$) is
said to be maximal (minimal). In addition, a set
$\{x_1,x_2,\cdots,x_k\}$ of points 
in P is said to be a chain if $x_{j+1}$ covers $x_j$ ($j=1,\cdots,
k-1$). A chain is 
maximal if $x_1$ and $x_k$ are respectively a minimal and a maximal point.

It is easy to read the topology of $P(M)$ once the partial order is given.
It is not difficult to check that
\[ O_x = \{ y \in P(M) \; : \; y \preceq x \} \; .\]
Indeed, one can even prove a stronger result \cite{So,Al},
namely that any finite set $P$ on
which a partial order $\preceq$ is defined can be made into a finite $T_0$
topological space by declaring that the smallest open neighborhood $O_x$
containing $x$ is given exactly by the above set. 

Throughout this article, we will use `finite poset' and `finite T$_0$ space'
interchangeably.

It is convenient to graphically represent a poset by a diagram, the Hasse
diagram, constructed by arranging its points at different levels and connecting
them using the following rules \cite{So,St}:\\
1. if $x\prec y$, then $x$ is at a lower level than $y$;\\
2. if $x\prec y$ and there is no $z$ such that $x\prec z\prec y$, then
$x$ is at a level immediately below $y$ and these two points are connected by
a line called a link.

Let us consider $P_4(S^1)$ again. The partial order reads
\be
x_1\preceq x_2~,~~ x_1\preceq x_4~,~~ x_3\preceq x_2~,~~
x_3\preceq x_4~,~~ \label{2.5}
\ee
where we have omitted writing the relations $x_j\preceq x_j$.
The corresponding Hasse diagram is shown in Figure \ref{fi:hassecircle}.

\begin{figure}[htb]
\begin{center}

\begin{picture}(300,120)(-150,-55)

\put(-40,40){\circle*{4}}
\put(40,40){\circle*{4}}
\put(-40,-40){\circle*{4}}
\put(40,-40){\circle*{4}}

\put(-40,40){\line(0,-1){80}}
\put(40,40){\line(0,-1){80}}
\put(-40,40){\line(1,-1){80}}
\put(40,40){\line(-1,-1){80}}

\put(-55,39){$x_4$}
\put(45,39){$x_2$}
\put(-55,-43){$x_1$}
\put(45,-43){$x_3$}

\end{picture}

\end{center}
\caption{\label{fi:hassecircle}}
\centerline{{\footnotesize  The
Hasse diagram for the circle poset $P_4 (S^1 )$.}}
\vskip1cm
\end{figure}

In the language of partially ordered sets, the smallest open set $O_x$
containing a point $x \in P(M)$ consists of all $y$ preceding
$x$: $O_x = \{y \in P(M) \; : \; y \preceq x \}$. In the Hasse diagram, it
consists of $x$ and all points we encounter as we travel along links
from $x$ to the bottom. In Figure \ref{fi:hassecircle}, this rule gives
$\{x_1,x_2,x_3\}$ as the smallest open set containing $x_2$, just as in
(\ref{2.4}).

As one example of a three-level poset, consider the Hasse diagram of Figure
\ref{fi:hassesphere} for a finite approximation $P_6(S^2)$ of the
two-dimensional 
sphere $S^2$ derived in \cite{So}. Its open sets are generated by
$$
\{x_1\}~,  ~ ~ \{x_3\}~,  ~ ~ \{x_1,x_2,x_3\}~,
~ ~ \{x_1,x_4,x_3\}~,
$$
\be
\{x_1,x_2,x_3,x_4,x_5\}~, ~ ~ \{x_1,x_2,x_3,x_4,x_6\}~,
~ ~ \label{2.6}
\ee
by taking unions and intersections.

\begin{figure}[htb]
\begin{center}

\begin{picture}(300,230)(-150,-140)

\put(-40,40){\circle*{4}}
\put(40,40){\circle*{4}}
\put(-40,-40){\circle*{4}}
\put(40,-40){\circle*{4}}
\put(-40,-120){\circle*{4}}
\put(40,-120){\circle*{4}}

\put(-40,40){\line(0,-1){80}}
\put(40,40){\line(0,-1){80}}
\put(-40,40){\line(1,-1){80}}
\put(40,40){\line(-1,-1){80}}
\put(-40,-40){\line(0,-1){80}}
\put(40,-40){\line(0,-1){80}}
\put(-40,-40){\line(1,-1){80}}
\put(40,-40){\line(-1,-1){80}}

\put(-55,39){$x_6$}
\put(45,39){$x_5$}
\put(-55,-43){$x_4$}
\put(45,-43){$x_2$}
\put(-55,-123){$x_1$}
\put(45,-123){$x_3$}

\end{picture}

\end{center}
\caption{\label{fi:hassesphere}}
\centerline{{\footnotesize  The
Hasse diagram for the sphere poset $P_6 (S^2)$.}}
\vskip1cm
\end{figure}

One of the most remarkable properties of a poset is its ability to accurately
reproduce the homology and the homotopy groups of the Hausdorff
topological space 
it approximates. For example, as for $S^1$, the fundamental group of $P_N(S^1)$
is ${\bf Z}$ whenever $N \geq 4$  \cite{So}. Similarly, as for $S^2$,
$\p_1(P_6(S^2)) =\{0\}$ and  $\p_2(P_6(S^2)) = {\bf Z}$. 
This has been widely discussed in our previous work \cite{BBET,lisbon},
where we argued that global topological information relevant for quantum
physics can be captured by such discrete approximations.
Furthermore, the topological space being
approximated can be recovered by considering a sequence of finer and finer
coverings, the appropriated framework being that of projective systems of
topological spaces. We refer to \cite{So,pangs} for details.

In this article we are however mostly
concerned with the algebraic properties of a poset, i.e. with the fact that any
finite poset can be regarded as the structure space of a
$C^*$-algebra. This will be 
extensively discussed and proved in the following sections, but let us first
illustrate a simple example.
Consider the \cstar :
\be
\ca= \{ \l_1 \unit_1 + \l_2 \unit_2 +
k_{12}\; : \; \l_j \in {\bf C}, k_{12}\in \ck_{12}\} \label{exa}
\ee
acting on the direct sum of two Hilbert spaces $\ch = \ch_1\bigoplus \ch_2$ and
generated by multiples of the identity $\unit_1$ on $\ch_1$, multiples of the
identity $\unit_2$ on $\ch_2$ and compact operators $\ck_{12}$ on the whole
Hilbert space $\ch$. This algebra admits only three classes of irreducible
representations, two finite dimensional ones and an infinite dimensional one:\\
1. $\p_1 : \l_1 \unit_1 + \l_2 \unit_2 + k_{12} \mapsto \l_1$,\\ 
2. $\p_2 : \l_1 \unit_1 + \l_2 \unit_2 + k_{12} \mapsto \l_2$,\\ 
3. $\r : \l_1 \unit_1 + \l_2 \unit_2 + k_{12}\mapsto \l_1 \unit_1 + \l_2
\unit_2 + k_{12}.$\\

\noindent
Hence the corresponding structure space consists of only three points
$p_1=ker\p_1,p_2=ker\p_2,q=ker\r$, corresponding respectively to the three
representations given above. This space has to be given the Jacobson
topology as 
explained in the previous section. This is easily done if one notices that, the
space being finite, this amounts to give a partial order relation on the set
$\{p_1,p_2,q\}$ \cite{FD}. Indeed one can show that, on any finite
structure space 
$Prim \ca$ of a \cstar $\ca$, the Jacobson topology is equivalent to the
following partial order relation:
\be
p_j \prec p_k \; \Leftrightarrow \; ker\p_j \subset ker\p_k\; ,
\ee
where $p_j$ is the point in $Prim\ca$ corresponding to the IRR $\p_j$ of $\ca$.
Thus in our example, since $ker\r  \subset ker\p_1$ and $ker\r
\subset ker\p_2$, the set $\{p_1,p_2,q\}$ is equipped with the order 
relations 
$q \prec p_1, q \prec p_2$ and therefore corresponds to the poset of Figure
\ref{fi:vposet}, which will be referred to as the $\bigvee$ poset from now on.

\begin{figure}[htb]
\begin{center}

\begin{picture}(320,150)(-60,20)

\put(100,50){\circle*{4}}
\put(50,100){\circle*{4}}
\put(150,100){\circle*{4}}

\put(100,50){\line(-1,1){50}}
\put(100,50){\line(1,1){50}}

\put(98,37){$q$}
\put(45,105){$p_1$}
\put(145,105){$p_2$}

\end{picture}

\end{center}
\caption{\label{fi:vposet}}
\centerline{{\footnotesize  The $\bigvee$ poset, structure space of
$\ca= {\bf C} \unit_1 + {\bf C} \unit_2 +
\ck_{12}$.}}
\vskip1cm
\end{figure}

\section{AF algebras}\label{se:af}

\subsection{Bratteli diagrams}\label{sub:dia}

A \cstar $\ca$ is said to be {\it{approximately finite dimensional}} (AF)
\cite{Br,Go} if there exists an increasing sequence 
\be
\ca_0 ~{\buildrel I_0 \over \hookrightarrow}~
\ca_1 ~{\buildrel I_1 \over \hookrightarrow}~ \ca_2
      ~{\buildrel I_2 \over \hookrightarrow}~ \ca_3
      ~{\buildrel I_2 \over \hookrightarrow}~ \cdots
      ~{\buildrel I_{n-1} \over \hookrightarrow}~ \ca_n
      ~{\buildrel I_n \over \hookrightarrow} \cdots
\label{af0}
\ee
of finite dimensional subalgebras of $\ca$, such that $\ca$ is the norm
closure of $\cup_n \ca_n$. Here the maps $I_n$ are injective
$^*$-homomorphisms. In other words, $\ca$ is the direct limit in the
category of \cstars with morphisms given by *-algebras maps (not isometries)
of the sequence $(\ca_n)_{n \in \IN}$.
As a set, $\bigcup_n \ca_n$ is made of coherent
sequences,
\be
\bigcup_n \ca_n = \{ a=(a_n)_{n \in \IN}~, a_n \in \ca_n ~|~ \exists  N_0 ~:
~a_{n+1} =  I_n(a_n)~, \forall ~n>N_0 \}.
\ee
Now the sequence $(||a_n||_{\ca_n})_{n \in \IN}$ is
eventually decreasing, since $||a_{n+1}|| \leq ||a_n||$ (the maps $I_n$ are
norm decreasing) and therefore convergent. One writes for the norm
\be
||(a_n)_{\IN}|| = \lim_{n \rt \infty} ||a_n||_{\ca_n}~. \label{norm1}
\ee
Since the maps $I_n$ are injective, the expression (\ref{norm1}) gives
directly a true norm and not simply a seminorm and there is no need to
quotient out the zero norm elements.

Each subalgebra $\ca_n$, being a finite dimensional $C^*$-algebra, is a
matrix algebra and therefore can be written as
$\ca_n = \bigoplus_{k=1}^{N_n} \IM^{(n)}(d_k,\IC)$ where
$\IM^{(n)}(d_k,\IC)$ is the algebra of $d_k \times d_k$ matrices
with complex coefficients. Given any two such matrix algebras $\ca_1
= \bigoplus_{j=1}^{N_1} \IM^{(1)}(d_j,\IC)$ and
$\ca_2=\bigoplus_{k=1}^{N_2} \IM^{(2)}(d_k,\IC)$ with $\ca_1 \hookrightarrow
\ca_2$, one can always choose suitable bases in $\ca_1$ and
$\ca_2$ such that $\ca_1$ is identified with a subalgebra of
$\ca_2$ of the following form \cite{Br}:
\be
\ca_1 \simeq \bigoplus_{k=1}^{N_2} \left( \bigoplus_{j=1}^{N_1} N_{kj}
\IM^{(1)}(d_j,\IC) \right)  \; .  \label{af}
\ee
Here, for any nonnegative integers $p$ and $q$, the symbol $p\,
\IM(q,\IC)$ stands 
for $\IM(q,\IC) \otimes \unit_p$. In (\ref{af}), the coefficients
$N_{kj}$ represent 
the multiplicity of the  partial embedding of $\IM^{(1)}(d_j,\IC)$ in
$\IM^{(2)}(d_k,\IC)$  and satisfy the condition 
\be
\sum_{j=1}^{N_1} N_{kj} d_j = d_k  \; . \label{dim}
\ee

A useful way to
represent the algebras $\ca_1$, $\ca_2$ and the embedding $\ca_1
\hookrightarrow \ca_2$ is by means of a diagram, the {\it Bratteli diagram}
\cite{Br}, which can be constructed out of the dimensions, 
$d_j ~(j=1,\ldots,N_1)$ and $d_k ~(k=1,\ldots,N_2)$, of the diagonal blocks of
the two algebras and the numbers $N_{kj}$ that describe the partial
embeddings. To construct the diagram, we draw two horizontal rows of vertices,
the top (bottom) one representing $\ca_1$ ($\ca_2$) and consisting of $N_1$
($N_2$) vertices, labeled by the corresponding dimensions
$d_1, \ldots, d_{N_1}$ ($d_1,\ldots,d_{N_2}$). Then for each
$j=1,\ldots,N_1$ and  
$k=1,\ldots,N_2$, we draw $N_{kj}$ edges between $d_j$ and $d_k$. We will also
write  $d^{(1)}_j \searrow^{N_{kj}} d^{(2)}_k$ to denote the fact that
$\IM(d^{(1)}_j,\IC)$ is embedded in $\IM(d^{(2)}_k,\IC)$ with multiplicity
$N_{kj}$.
By repeating the procedure at each level, we obtain a semi-infinite
diagram denoted by $\cd(\ca)$ which completely defines $\ca$ up to
isomorphisms. 
Notice that the diagram $\cd(\ca)$ depends not only on $\ca$ but also on the
particular sequence $\{\ca_n\}_{n \in \IN}$ which generates
$\ca$. However, it is possible to show \cite{Br} that all diagrams
corresponding 
to AF algebras which are isomorphic to $\ca$ can be obtain from the chosen 
$\cd(\ca)$ by means of an algorithm.

As an example of an AF algebra, let us consider the subalgebra $\ca$ of the
algebra ${\cal B}(\ch )$ of bounded operators on $\ch = \ch_1 \oplus
\ch_2$ given
in (\ref{exa}). This \cstar algebra can be obtained as the direct limit of the
following sequence of finite dimensional algebras:\\
\vskip.3cm
$\begin{array}{l}
\ca_0 = \IM(1,\IC) \\
\ca_1 = \IM(1,\IC) \oplus \IM(1,\IC) \\
\ca_2 = \IM(1,\IC) \oplus \IM(2,\IC) \oplus \IM(1,\IC)\\
\vdots \\
\ca_n = \IM(1,\IC) \oplus \IM(2n-2,\IC) \oplus \IM(1,\IC)\\
\vdots
\end{array}$\\
\vskip.3cm
\noindent
where, for $n \geq 1$, $\ca_n$ is embedded in $\ca_{n+1}$ as the 
subalgebra $\IM(1,\IC) \oplus
\left[ \IM(1,\IC) \oplus \IM(2n-2,\IC) \oplus \IM(1,\IC) \right] \oplus
\IM(1,\IC)$:
\be a_{n} = \left[
\begin{array}{ccc}
\l_1 & 0                  & 0    \\
0    & m_{2n-2 \times 2n-2} & 0    \\
0    & 0                  & \l_2
\end{array}
\right]~~ \hookrightarrow ~
\left[
\begin{array}{ccccc}
\l_1 & 0    & 0                  & 0       & 0      \\
0    & \l_1 & 0                  & 0       & 0      \\
0    & 0    & m_{2n-2 \times 2n-2} & 0       & 0      \\
0    & 0    & 0                  & \l_2    & 0       \\
0    & 0    & 0                  & 0       & \l_2
\end{array}
\right]~.
\label{vee1}
\ee

It is therefore described by the diagram of Figure \ref{fi:valg}.

\begin{figure}[htb]
\begin{center}
\begin{picture}(120,200)(0,20)

\put(30,150){\circle*{4}}
\put(30,120){\circle*{4}}
\put(30,90){\circle*{4}}
\put(30,60){\circle*{4}}

\put(60,180){\circle*{4}}
\put(60,120){\circle*{4}}
\put(60,90){\circle*{4}}
\put(60,60){\circle*{4}}

\put(90,150){\circle*{4}}
\put(90,120){\circle*{4}}
\put(90,90){\circle*{4}}
\put(90,60){\circle*{4}}

\put(30,150){\line(0,-1){30}}
\put(30,120){\line(0,-1){30}}
\put(30,90){\line(0,-1){30}}
\put(60,120){\line(0,-1){30}}
\put(60,90){\line(0,-1){30}}
\put(90,150){\line(0,-1){30}}
\put(90,120){\line(0,-1){30}}
\put(90,90){\line(0,-1){30}}
\put(30,60){\line(0,-1){10}}
\put(60,60){\line(0,-1){10}}
\put(90,60){\line(0,-1){10}}

\put(60,180){\line(1,-1){30}}
\put(30,150){\line(1,-1){30}}
\put(30,120){\line(1,-1){30}}
\put(30,90){\line(1,-1){30}}
\put(60,180){\line(-1,-1){30}}
\put(90,150){\line(-1,-1){30}}
\put(90,120){\line(-1,-1){30}}
\put(90,90){\line(-1,-1){30}}
\put(30,60){\line(1,-1){10}}
\put(90,60){\line(-1,-1){10}}

\put(58,183){1}

\put(22,148){1}
\put(22,118){1}
\put(22,88){1}
\put(22,58){1}

\put(93,148){1}
\put(93,118){1}
\put(93,88){1}
\put(93,58){1}

\put(62,113){2}
\put(62,83){4}
\put(62,53){6}

\put(45,30){$\vdots$}
\put(75,30){$\vdots$}

\end{picture}
\end{center}
\caption{\label{fi:valg}}
\centerline{{\footnotesize The Bratteli diagram corresponding to $\ca = \IC
\unit_1 +\IC \unit_2 +\ck_{12}$.}}
\vskip1cm
\end{figure}

As a second example, consider the \cstar of the Penrose tiling. This is an
example of an AF algebra which is not postliminal,  since this algebra
admits an infinite number of nonequivalent representations all with
the same kernel. At each 
level, the finite dimensional algebra is given by \cite{Co}
\be
\ca_n = \IM(d_n,\IC) \oplus \IM(d'_n,\IC)~, ~~ n \geq 1 ~, \label{pt}
\ee
with inclusion $\ca_n \hookrightarrow \ca_{n+1}$:
\be
\left[
\begin{array}{cc}
A & 0 \\
0 & B
\end{array}
\right] \hookrightarrow
\left[
\begin{array}{ccc}
A & 0 & 0 \\
0 & B & 0 \\
0 & 0 & A
\end{array}
\right] ~; \label{ptinc}
~~A \in \IM(d_n,\IC)~, ~~B \in \IM(d'_n,\IC)~,
\ee
so that $d_{n+1}= d_n + d_n'$ and $d_{n+1}'=d_n$.
The corresponding Bratteli diagram is shown in Figure \ref{fi:penbra}.

\begin{figure}[htb]
\begin{center}
\begin{picture}(150,190)(-10,-10)

\put(60,150){\circle*{4}}
\put(30,120){\circle*{4}}
\put(30,90){\circle*{4}}
\put(30,60){\circle*{4}}
\put(30,30){\circle*{4}}
\put(90,120){\circle*{4}}
\put(90,90){\circle*{4}}
\put(90,60){\circle*{4}}
\put(90,30){\circle*{4}}

\put(60,150){\line(-1,-1){30}}
\put(60,150){\line(1,-1){30}}
\put(30,120){\line(0,-1){30}}
\put(30,90){\line(0,-1){30}}
\put(30,60){\line(0,-1){30}}
\put(30,120){\line(2,-1){60}}
\put(30,90){\line(2,-1){60}}
\put(30,60){\line(2,-1){60}}
\put(90,120){\line(-2,-1){60}}
\put(90,90){\line(-2,-1){60}}
\put(90,60){\line(-2,-1){60}}
\put(30,30){\line(2,-1){20}}
\put(30,30){\line(0,-1){10}}
\put(90,30){\line(-2,-1){20}}

\put(60,0){$\vdots$}
\put(57,154){$1$}
\put(22,118){$1$}
\put(22,88){$2$}
\put(22,58){$3$}
\put(22,28){$5$}
\put(94,118){$1$}
\put(94,88){$1$}
\put(94,58){$2$}
\put(94,28){$3$}

\end{picture}
\end{center}
\caption{\label{fi:penbra}}
\centerline{{\footnotesize The Bratteli diagram of the Penrose tiling.}}
\vskip1cm
\end{figure}

To conclude this section we remark that an AF algebra is commutative
if and only 
if all its factors $\IM^{(n)}(d_k,\IC)$ are one dimensional, i.e. they are just
$\IC$. Thus the corresponding diagram has the property that for each
$\IM^{(n)}(d_k,\IC)$ with $n \geq 1$ there is exactly one
$\IM^{(n-1)}(d_j,\IC)$ and 
$\IM^{(n-1)}(d_j,\IC) \searrow^{p_{kj}} \IM^{(n)}(d_k,\IC)$ with
$p_{kj}=1$. An interesting example is given in Figure \ref{fi:comalg}, which
corresponds to the AF \cstar of continuous functions on the Cantor set
\cite{Ef}. 

\begin{figure}[htb]
\begin{center}

\begin{picture}(400,200)(-200,-140)
\put(0,10){\circle*{4}}
\put(-20,-90){\circle*{4}}
\put(-40,-60){\circle*{4}}
\put(-60,-90){\circle*{4}}
\put(-80,-30){\circle*{4}}
\put(-100,-90){\circle*{4}}
\put(-120,-60){\circle*{4}}
\put(-140,-90){\circle*{4}}
\put(20,-90){\circle*{4}}
\put(40,-60){\circle*{4}}
\put(60,-90){\circle*{4}}
\put(80,-30){\circle*{4}}
\put(100,-90){\circle*{4}}
\put(120,-60){\circle*{4}}
\put(140,-90){\circle*{4}}
\put(-142,-110){$\vdots$}
\put(-102,-110){$\vdots$}
\put(-62,-110){$\vdots$}
\put(-22,-110){$\vdots$}
\put(138,-110){$\vdots$}
\put(98,-110){$\vdots$}
\put(58,-110){$\vdots$}
\put(18,-110){$\vdots$}

\put(0,10){\line(-2,-1){80}}
\put(0,10){\line(2,-1){80}}

\put(80,-30){\line(-4,-3){40}}
\put(80,-30){\line(4,-3){40}}
\put(-80,-30){\line(-4,-3){40}}
\put(-80,-30){\line(4,-3){40}}

\put(120,-60){\line(2,-3){20}}
\put(120,-60){\line(-2,-3){20}}
\put(-120,-60){\line(2,-3){20}}
\put(-120,-60){\line(-2,-3){20}}

\put(40,-60){\line(-2,-3){20}}
\put(40,-60){\line(2,-3){20}}
\put(-40,-60){\line(-2,-3){20}}
\put(-40,-60){\line(2,-3){20}}

\put(140,-90){\line(1,-1){10}}
\put(140,-90){\line(-1,-1){10}}
\put(-140,-90){\line(1,-1){10}}
\put(-140,-90){\line(-1,-1){10}}

\put(100,-90){\line(1,-1){10}}
\put(100,-90){\line(-1,-1){10}}
\put(-100,-90){\line(1,-1){10}}
\put(-100,-90){\line(-1,-1){10}}

\put(60,-90){\line(1,-1){10}}
\put(60,-90){\line(-1,-1){10}}
\put(-60,-90){\line(1,-1){10}}
\put(-60,-90){\line(-1,-1){10}}

\put(20,-90){\line(1,-1){10}}
\put(20,-90){\line(-1,-1){10}}
\put(-20,-90){\line(1,-1){10}}
\put(-20,-90){\line(-1,-1){10}}

\put(-1,14){1}
\put(12,-92){1}
\put(-16,-92){1}
\put(32,-62){1}
\put(-37,-62){1}
\put(52,-92){1}
\put(-55,-92){1}
\put(82,-30){1}
\put(-88,-30){1}
\put(92,-92){1}
\put(-97,-92){1}
\put(123,-62){1}
\put(-128,-62){1}
\put(143,-92){1}
\put(-148,-92){1}

\end{picture}
\end{center}
\caption{\label{fi:comalg}}
\centerline{{\footnotesize The Bratteli diagram corresponding to the
AF \cstar}} \centerline{{\footnotesize of continuous functions on the
Cantor set.}} 
\vskip1cm
\end{figure}

\subsection{From Bratteli Diagrams to Posets}\label{sub:brapo}

The Bratteli diagram $\cd(\ca)$ of an AF algebra $\ca$ is useful not only
because it gives the finite approximations of the algebra explicitly,
but also because it is possible to read the ideals and the primitive
ideals of the algebra (hence the topological properties of $Prim\ca$) out of it
very easily. Indeed one can show that the following proposition holds
\cite{Br}:
\bprop \label{th:braideal}
\begin{enumerate}
\item There is a one-to-one correspondence between the proper ideals
$\ci$ of $\ca$ and the subsets $\Lambda = \Lambda_{\ci}$ of the
Bratteli diagram 
satisfying the following two properties:\\
$i$) if $\IM^{(n)}(d_k,\IC)\in \Lambda$ and $\IM^{(n)}(d_k,\IC) \searrow
\IM^{(n+1)}(d_j,\IC)$ then necessarily $\IM^{(n+1)}(d_j,\IC)$ belongs to
$\Lambda$;\\
$ii$) if all factors $\IM^{(n+1)}(d_j,\IC)$
($j=\{1,2,\cdots,N_{n+1}\}$), for which 
$\IM^{(n)}(d_k,\IC) \searrow \IM^{(n+1)}(d_j,\IC)$, belong to $\Lambda$, then 
$\IM^{(n)}(d_k,\IC) \in \Lambda$.\\
\item A proper ideal $\ci$ of $\ca$ is primitive if and only if the associated
subdiagram $\Lambda_{\ci}$ satisfies:\\
$iii$) $\forall n$ there exists an $\IM^{(m)}(d_j,\IC)$, with $m > n$, not
belonging to $\Lambda_{\ci}$ such that, for all $k\in\{1,2,\cdots,N_n\}$ with
$\IM^{(n)}(d_k,\IC)$ not in $\Lambda_{\ci}$, one can find a sequence
$\IM^{(n)}(d_k,\IC) \searrow \IM^{(n+1)}(d_{\a},\IC) \searrow
\IM^{(n+2)}(d_{\b},\IC) \searrow \cdots \searrow \IM^{(m)}(d_j,\IC)$.
\end{enumerate}
\eprop

For example, consider the diagram of Figure \ref{fi:valg}, representing the
AF \cstar $\ca = \IC \unit_1 +\IC \unit_2 +\ck_{12}$ already discussed in
section \ref{sub:fta}. This algebra contains only three
nontrivial ideals, whose diagrams are represented in Figure
\mbox{\ref{fi:videal}(a,b,c)}. In this pictures the points belonging
to the ideals 
are marked with a ``$\clubsuit$". It is not difficult to check that
only $\ci_1$ and 
$\ci_2$ are primitive ideals, since $\ci_3$ does not satisfy property $(iii)$
above.\\

\begin{figure}[htb]
\begin{center}
\begin{picture}(360,230)(0,0)

\put(30,150){\circle*{4}}
\put(30,120){\circle*{4}}
\put(30,90){\circle*{4}}
\put(30,60){\circle*{4}}

\put(60,180){\circle*{4}}
\put(56,117){$\clubsuit$}
\put(56,87){$\clubsuit$}
\put(56,57){$\clubsuit$}

\put(86,147){$\clubsuit$}
\put(86,117){$\clubsuit$}
\put(86,87){$\clubsuit$}
\put(86,57){$\clubsuit$}

\put(30,150){\line(0,-1){30}}
\put(30,120){\line(0,-1){30}}
\put(30,90){\line(0,-1){30}}
\put(60,120){\line(0,-1){30}}
\put(60,90){\line(0,-1){30}}
\put(90,150){\line(0,-1){30}}
\put(90,120){\line(0,-1){30}}
\put(90,90){\line(0,-1){30}}
\put(30,60){\line(0,-1){10}}
\put(60,60){\line(0,-1){10}}
\put(90,60){\line(0,-1){10}}

\put(60,180){\line(1,-1){30}}
\put(30,150){\line(1,-1){30}}
\put(30,120){\line(1,-1){30}}
\put(30,90){\line(1,-1){30}}
\put(60,180){\line(-1,-1){30}}
\put(90,150){\line(-1,-1){30}}
\put(90,120){\line(-1,-1){30}}
\put(90,90){\line(-1,-1){30}}
\put(30,60){\line(1,-1){10}}
\put(90,60){\line(-1,-1){10}}

\put(146,147){$\clubsuit$}
\put(146,117){$\clubsuit$}
\put(146,87){$\clubsuit$}
\put(146,57){$\clubsuit$}

\put(180,180){\circle*{4}}
\put(176,117){$\clubsuit$}
\put(176,87){$\clubsuit$}
\put(176,57){$\clubsuit$}

\put(210,150){\circle*{4}}
\put(210,120){\circle*{4}}
\put(210,90){\circle*{4}}
\put(210,60){\circle*{4}}

\put(150,150){\line(0,-1){30}}
\put(150,120){\line(0,-1){30}}
\put(150,90){\line(0,-1){30}}
\put(180,120){\line(0,-1){30}}
\put(180,90){\line(0,-1){30}}
\put(210,150){\line(0,-1){30}}
\put(210,120){\line(0,-1){30}}
\put(210,90){\line(0,-1){30}}
\put(150,60){\line(0,-1){10}}
\put(180,60){\line(0,-1){10}}
\put(210,60){\line(0,-1){10}}

\put(180,180){\line(1,-1){30}}
\put(150,150){\line(1,-1){30}}
\put(150,120){\line(1,-1){30}}
\put(150,90){\line(1,-1){30}}
\put(180,180){\line(-1,-1){30}}
\put(210,150){\line(-1,-1){30}}
\put(210,120){\line(-1,-1){30}}
\put(210,90){\line(-1,-1){30}}
\put(150,60){\line(1,-1){10}}
\put(210,60){\line(-1,-1){10}}

\put(270,150){\circle*{4}}
\put(270,120){\circle*{4}}
\put(270,90){\circle*{4}}
\put(270,60){\circle*{4}}

\put(300,180){\circle*{4}}
\put(296,117){$\clubsuit$}
\put(296,87){$\clubsuit$}
\put(296,57){$\clubsuit$}

\put(330,150){\circle*{4}}
\put(330,120){\circle*{4}}
\put(330,90){\circle*{4}}
\put(330,60){\circle*{4}}

\put(270,150){\line(0,-1){30}}
\put(270,120){\line(0,-1){30}}
\put(270,90){\line(0,-1){30}}
\put(300,120){\line(0,-1){30}}
\put(300,90){\line(0,-1){30}}
\put(330,150){\line(0,-1){30}}
\put(330,120){\line(0,-1){30}}
\put(330,90){\line(0,-1){30}}
\put(270,60){\line(0,-1){10}}
\put(300,60){\line(0,-1){10}}
\put(330,60){\line(0,-1){10}}

\put(300,180){\line(1,-1){30}}
\put(270,150){\line(1,-1){30}}
\put(270,120){\line(1,-1){30}}
\put(270,90){\line(1,-1){30}}
\put(300,180){\line(-1,-1){30}}
\put(330,150){\line(-1,-1){30}}
\put(330,120){\line(-1,-1){30}}
\put(330,90){\line(-1,-1){30}}
\put(270,60){\line(1,-1){10}}
\put(330,60){\line(-1,-1){10}}

\put(52,15){$(a)$}
\put(172,15){$(b)$}
\put(292,15){$(c)$}
\put(75,170){$\ci_1$}
\put(195,170){$\ci_2$}
\put(315,170){$\ci_3$}

\put(45,30){$\vdots$}
\put(75,30){$\vdots$}
\put(165,30){$\vdots$}
\put(195,30){$\vdots$}
\put(285,30){$\vdots$}
\put(315,30){$\vdots$}

\end{picture}
\end{center}
\caption{\label{fi:videal}}
\centerline{{\footnotesize The representation of the ideals of
$\ca = \IC \unit_1 +\IC \unit_2 +\ck_{12}$}}
\centerline{{\footnotesize in the corresponding Bratteli
diagram.}}
\vskip1cm
\end{figure}

We remark the following here:\\
1) The whole $\ca$ is an ideal, which by definition is not primitive since the
trivial representation $\ca \rightarrow 0$ is not irreducible.\\
2) The set $\{0\} \subset \ca$ is an ideal, which is primitive if and only if
$\ca$ has one irreducible faithful representation. This can also be understood
from the Bratteli diagram in the following way. The set $\{0\}$ is not a
subdiagram of $\cd(\ca)$, being represented by the element $0$ of the matrix
algebra of each finite level, so that there is at least one element $a\in \ca$
not belonging to the ideal $\{0\}$ at any level. Thus to check if $\{0\}$ is
primitive, i.e. to check property ($iii$) above,
we have to examine whether {\it all} the points at a
given level, say $n$, can be connected to a {\it single} point at a
level $m>n$.
For example this is the case for the diagram of Figure \ref{fi:valg} and not
for that of Figure \ref{fi:comalg}.  \\

Proposition \ref{th:braideal} above allows us to understand the topological
properties of $Prim\ca$ at once. This becomes particularly simple if the
algebra admits only  a finite number of nonequivalent irreducible
representations. In this case $Prim\ca$ is a $T_0$ topological space with only
a finite number of points, hence a finite poset $P$. To reconstruct the latter
we just need to draw the Bratteli diagram $\cd(\ca)$ and find the subdiagrams
that, according to properties ($i,ii,iii$), correspond to primitive
ideals. Then 
$P$ has so many points as the number of primitive ideals and the partial order
relation in $P$ that determines the $T_0$ topology is simply given by the
inclusion relations that exist among the primitive ideals.

As an example consider again Figure \ref{fi:valg}. We have seen that the
corresponding AF algebra has only three primitive ideals: the $\{0\}$ ideal and
the ideals $\ci_1,\ci_2$ represented in Figure \ref{fi:videal}(a),(b). Clearly
$\{0\}\subset\ci_1,\ci_2$ so that $Prim\ca$ is the $\bigvee$ poset of Figure
\ref{fi:vposet}.

Figure \ref{fi:comalg} leads to another interesting topological
space. As we have 
mentioned, such a diagram corresponds to a commutative AF algebra $\cc$ and
hence to a Hausdorff $Prim\cc$, which is homeomorphic to the Cantor set.

\subsection{From Posets to Bratteli Diagrams}\label{sub:pobra}

In the preceding subsection we have described the properties
of the Bratteli diagram $\cd(\ca)$ of an AF algebra $\ca$ and in particular we
have seen how, out of it, it is possible to read the primitive ideal space of
$\ca$, in particular when the latter is a finite poset. In the
following we will 
see that, under some rather mild hypotheses which are always verified in the
cases of interest to us, it is possible to reverse the construction and thus
build the AF algebra that corresponds to a given (finite)  $T_0$ topological
space.

Such a reconstruction rests on the following theorem of Bratteli
\cite{Br}, which 
specifies a class of topological spaces which are the primitive ideal
spaces of AF algebras:
\bprop \label{th:bratop}
A topological space $Y$ is the primitive ideal space $Prim\ca$
of an AF algebra $\ca$ if it has the following properties:\\
i) $Y$ is $T_0$;\\
ii) $Y$ contains at most a countable number of closed sets;\\
iii) if $\{F_n\}_{n\in \L}$, $\L$ being any direct set, is a
decreasing sequence of 
closed subsets of $Y$, then $\cap_n F_n$ is an element in
$\{F_n\}_{n\in \L}$;\\ 
iv) if $F\subset Y$ is a closed set which is not the union of two proper closed
subsets, then $F$ is the closure of a one-point set.
\eprop

It is not difficult to check that all the above conditions hold true
if $Y$ is a 
$T_0$ topological space with a finite number of points, so that we have the
corollary:\\

\noindent
{\bf Corollary.}\\
 {\it A finite poset $P$ is the primitive ideal space $Prim\ca$ for
some AF algebra $\ca$.}\\

Here we will not report the proof of proposition \ref{th:bratop}, which can be
found in \cite{Br}. However, starting from the techniques used in such
a proof, we 
want to show how one can explicitly find an AF algebra $\ca$ whose primitive
ideal space is a given finite poset $P$. First we will give the general
construction and then discuss an example.

Let $\{K_1,K_2,K_3,\ldots\}$ be the collection of all closed sets in $P$,
where $K_1 = P$. To construct the $n$-th level of the Bratteli diagram
$\cd(\ca)$, we consider the subcollection
of closed sets $\ck_n\equiv \{K_1,K_2,\ldots ,K_n\}$
and denote with $\ck_n'$ the smallest collection of closed sets in $P$
that contains $\ck_n$ and is closed under union and intersection. The
collection 
$\ck_n$ determines a partition of the topological space $P$, 
by taking intersections and complements of the sets $K_j \in \ck_n$
($j=1,\ldots,n$). We denote with $Y(n,1),\; Y(n,2), \; \ldots ,\;Y(n,k_n)$ the
sets of such partition. Also, we write $F(n,j)$ for the smallest
closed set which 
contains $Y(n,j)$ and belongs to the subcollection $\ck_n'$. Then we
can construct a Bratteli diagram following the rules:
\begin{enumerate}
\item the $n$-th level of $\cd(\ca)$ has $k_n$ points, one for each set
$Y(n,j)$;
\item the point at the level $n$ of the diagram corresponding to
$Y(n,\a)$ is linked to the point at the level $n+1$ corresponding to
$Y(n+1,\b)$ 
if and only if $Y(n,\a)\cap F(n+1,\b) \neq \emptyset$. In this case, the
multiplicity of the embedding is always 1.
\end{enumerate}

To illustrate this construction, let us consider the $\bigvee$ poset of Figure
\ref{fi:vposet}: $P=\{p_1,p_2,q\}$. Now there are four closed sets:\\
\[ K_1=\{p_1,p_2,q\},\; K_2=\{p_1\},\; K_3=\{p_2\},\; K_4=\{p_1,p_2\}\; . \]
Thus it is not difficult to check that:
\begin{displaymath}
\begin{small}
\begin{array}{lllllll}
\ck_1=\{K_1\}             & & \ck_1'=\{K_1\}             & &
                                                       Y(2,1)=\{p_1,p_2,q\} 
                                                       &\subset& F(1,1)=K_1\\
 & & & & & & \\
\ck_2=\{K_1,K_2\}         & &\ck_2'=\{K_1,K_2\}          & & Y(2,1)=\{p_1\}
                                                       &\subset& F(2,1)=K_2\\
                          & &                            & & Y(2,2)=\{p_2,q\}
                                                       &\subset& F(2,2)=K_1\\
 & & & & & & \\
\ck_3=\{K_1,K_2,K_3\}     & & \ck_3'=\{K_1,K_2,K_3,K_4\} & & Y(3,1)=\{p_1\}
                                                       &\subset& F(3,1)=K_2 \\
                          & &                            & & Y(3,2)=\{q\}
                                                       &\subset& F(3,2)=K_1 \\
                          & &                            & & Y(3,3)=\{p_2\}
                                                       &\subset& F(3,3)=K_3 \\
 & & & & & & \\
\ck_4=\{K_1,K_2,K_3,K_4\} & & \ck_4'=\{K_1,K_2,K_3,K_4\} & & Y(4,1)=\{p_1\}
                                                       &\subset& F(4,1)=K_2 \\
                          & &                            & & Y(4,2)=\{q\}
                                                       &\subset& F(4,2)=K_1 \\
                          & &                            & & Y(4,3)=\{p_2\}
                                                       &\subset& F(4,3)=K_3 \\
 & & & \vdots & & &
\end{array}
\end{small}
\end{displaymath}

Notice that, since $P$ has only a finite number of points and hence a finite
number of closed sets, the partition of $P$ we have to consider at each level
$n$ repeats itself after a certain point ($n=3$ in this case). Figure
\ref{fi:vtop} shows the corresponding diagram, obtained through rules (1) and
(2) above. Recalling then that the first matrix algebra that gives
origin to an AF 
algebra is $\IC$ and using the fact that all the embeddings have
multiplicity one, we eventually obtain the sequence of finite dimensional
algebras shown by the Bratteli diagram of Figure \ref{fi:valg}. As we have said
previously,  such a diagram corresponds to the AF algebra  $\ca = \IC
\unit_1 +\IC \unit_1 + \ck_{12}$.

\begin{figure}[htb]
\begin{center}
\begin{picture}(120,220)(0,20)

\put(30,150){\circle*{4}}
\put(30,120){\circle*{4}}
\put(30,90){\circle*{4}}
\put(30,60){\circle*{4}}

\put(60,180){\circle*{4}}
\put(60,120){\circle*{4}}
\put(60,90){\circle*{4}}
\put(60,60){\circle*{4}}

\put(90,150){\circle*{4}}
\put(90,120){\circle*{4}}
\put(90,90){\circle*{4}}
\put(90,60){\circle*{4}}

\put(30,150){\line(0,-1){30}}
\put(30,120){\line(0,-1){30}}
\put(30,90){\line(0,-1){30}}
\put(60,120){\line(0,-1){30}}
\put(60,90){\line(0,-1){30}}
\put(90,150){\line(0,-1){30}}
\put(90,120){\line(0,-1){30}}
\put(90,90){\line(0,-1){30}}
\put(30,60){\line(0,-1){10}}
\put(60,60){\line(0,-1){10}}
\put(90,60){\line(0,-1){10}}

\put(60,180){\line(1,-1){30}}
\put(30,150){\line(1,-1){30}}
\put(30,120){\line(1,-1){30}}
\put(30,90){\line(1,-1){30}}
\put(60,180){\line(-1,-1){30}}
\put(90,150){\line(-1,-1){30}}
\put(90,120){\line(-1,-1){30}}
\put(90,90){\line(-1,-1){30}}
\put(30,60){\line(1,-1){10}}
\put(90,60){\line(-1,-1){10}}

\put(54,185){$Y_{11}$}

\put(14,148){$Y_{21}$}
\put(14,118){$Y_{31}$}
\put(14,88){$Y_{41}$}
\put(14,58){$Y_{51}$}

\put(93,148){$Y_{22}$}
\put(93,118){$Y_{33}$}
\put(93,88){$Y_{43}$}
\put(93,58){$Y_{53}$}

\put(63,113){$Y_{32}$}
\put(63,83){$Y_{42}$}
\put(63,53){$Y_{51}$}

\put(45,30){$\vdots$}
\put(75,30){$\vdots$}

\end{picture}
\end{center}
\caption{\label{fi:vtop}}
\centerline{{\footnotesize The construction of the Bratteli diagram
of the AF algebra }}
\centerline{{\footnotesize corresponding to the $\bigvee$ poset of Figure
\ref{fi:vposet}.}}
\vskip1cm
\end{figure}

It is a general fact that the Bratteli diagram describing any finite poset
``stabilizes", i.e. repeats itself, after a certain level $n_0$, when
the family 
$\ck_{n_0}$ of closed sets we choose is such that it determines a
partition of the 
poset which distinguishes each point of the poset itself. In
particular, this is 
the case if we choose $n_0$ in such a manner that $\ck_{n_0}$ contains
all closed 
sets. Then, each $Y(n_0,j)$ will contain a single point of the poset and
$F(n_0+1,j)$ will be the smallest closed set containing $Y(n_0,j)$. It is only
this stable part of the diagram which is relevant for the inductive limit and
hence for the determination of the AF algebra it represents. Indeed,
diagrams (or 
sequences of finite dimensional algebras) that differ only for a
finite numbers of 
initial levels give different finite approximations to the {\underline
{same}} AF 
algebra \cite{Br,Go}.

To conclude this section, we want to describe the AF algebras whose structure
spaces are the poset approximations of the circle, $P_4(S^1)$, and of
the sphere, 
$P_6(S^2)$.

As for  $P_4(S^1)$, given in Figure \ref{fi:hassecircle}, the Bratteli diagram
repeats itself for $n>n_0 =4$ and the stable partition is given by 
\be
\begin{array}{ll}
Y(n_0,1)=\{x_2 \}  ~~&~~ F(n_0+1,1)=\{x_2 \} \\
Y(n_0,2)=\{x_1 \}  ~~&~~ F(n_0+1,2)=\{x_1, x_2, x_4 \} \\
Y(n_0,3)=\{x_3 \}  ~~&~~ F(n_0+1,3)=\{x_2, x_3, x_4 \} \\
Y(n_0,4)=\{x_4 \}  ~~&~~ F(n_0+1,4)=\{x_4 \}~.
\end{array}
\ee
\begin{figure}[htb]
\begin{center}
\begin{picture}(120,150)(10,20)
\put(30,120){\circle*{4}}
\put(30,90){\circle*{4}}
\put(30,60){\circle*{4}}
\put(60,120){\circle*{4}}
\put(60,90){\circle*{4}}
\put(60,60){\circle*{4}}
\put(90,120){\circle*{4}}
\put(90,90){\circle*{4}}
\put(90,60){\circle*{4}}
\put(120,120){\circle*{4}}
\put(120,90){\circle*{4}}
\put(120,60){\circle*{4}}
\put(30,120){\line(0,-1){30}}
\put(30,90){\line(0,-1){30}}
\put(60,120){\line(0,-1){30}}
\put(60,90){\line(0,-1){30}}
\put(90,120){\line(0,-1){30}}
\put(90,90){\line(0,-1){30}}
\put(120,120){\line(0,-1){30}}
\put(120,90){\line(0,-1){30}}
\put(30,120){\line(1,-1){30}}
\put(30,90){\line(1,-1){30}}
\put(120,120){\line(-1,-1){30}}
\put(120,90){\line(-1,-1){30}}
\put(120,120){\line(-2,-1){60}}
\put(120,90){\line(-2,-1){60}}
\put(30,120){\line(2,-1){60}}
\put(30,90){\line(2,-1){60}}
\put(30,60){\line(0,-1){10}}
\put(30,60){\line(1,-1){10}}
\put(60,60){\line(0,-1){10}}
\put(90,60){\line(0,-1){10}}
\put(120,60){\line(-1,-1){10}}
\put(120,60){\line(0,-1){10}}
\put(120,60){\line(-2,-1){20}}
\put(30,120){\line(0,1){10}}
\put(60,120){\line(0,1){10}}
\put(60,120){\line(-1,1){10}}
\put(90,120){\line(0,1){10}}
\put(90,120){\line(1,1){10}}
\put(120,120){\line(0,1){10}}
\put(60,120){\line(2,1){20}}
\put(90,120){\line(-2,1){20}}
\put(75,140){$\vdots$}
\put(75,30){$\vdots$}
\end{picture}
\end{center}
\caption{\label{fi:cirbra}} 
\centerline{\footnotesize{The stable part of the Bratteli diagram
for the circle poset $P_4(S^1)$.}}
\vskip1cm
\end{figure}
The corresponding Bratteli diagram is in Figure \ref{fi:cirbra}.
The set $\{0\}$ is not an ideal. The limit algebra $\ca$ turns out to
be a subalgebra of bounded operators on the Hilbert space $\ch = \ch_1 \oplus
\cdots \oplus \ch_4$, with $\ch_i ,~ i = 1, \dots, 4$ infinite
dimensional Hilbert spaces:
\be
\ca = \IC \unit_{13} \oplus \IC \unit_{24}
\oplus \ck_{12} \oplus \ck_{34} \;.
\label{ciralg}
\ee
Here $\unit_{ij}$ and $\ck_{ij}$ denote the identity operator and the
algebra of 
compact operators on $\ch_i \oplus \ch_j$ respectively.

For the poset $P_6(S^2)$ for the two-dimensional sphere, given in
Figure \ref{fi:hassesphere}, $n_0=6$ and the stable partition is given by 
\be
\begin{array}{ll}
Y(n_0,1)=\{x_5 \}  ~~&~~ F(n_0+1,1)=\{x_5 \}\\
Y(n_0,2)=\{x_2 \}  ~~&~~ F(n_0+1,2)=\{x_2, x_5, x_6 \} \\
Y(n_0,3)=\{x_1 \}  ~~&~~ F(n_0+1,3)=\{x_1, x_2, x_4, x_5, x_6, \} \\
Y(n_0,4)=\{x_3 \}  ~~&~~ F(n_0+1,4)=\{x_2, x_3, x_4, x_5, x_6 \} \\
Y(n_0,5)=\{x_4 \}  ~~&~~ F(n_0+1,5)=\{x_4, x_5, x_6 \} \\
Y(n_0,6)=\{x_6 \}  ~~&~~ F(n_0+1,6)=\{x_6 \} ~.
\end{array}
\ee
The corresponding Bratteli diagram is in Figure \ref{fi:sphbra}.
\begin{figure}[htb]
\begin{center}
\begin{picture}(600,300)(-200,40)
\put(-175,100){\circle*{4}}
\put(-175,170){\circle*{4}}
\put(-175,240){\circle*{4}}
\put(-105,100){\circle*{4}}
\put(-105,170){\circle*{4}}
\put(-105,240){\circle*{4}}
\put(-35,100){\circle*{4}}
\put(-35,170){\circle*{4}}
\put(-35,240){\circle*{4}}
\put(175,100){\circle*{4}}
\put(175,170){\circle*{4}}
\put(175,240){\circle*{4}}
\put(105,100){\circle*{4}}
\put(105,170){\circle*{4}}
\put(105,240){\circle*{4}}
\put(35,100){\circle*{4}}
\put(35,170){\circle*{4}}
\put(35,240){\circle*{4}}
\put(-175,170){\line(0,-1){70}}
\put(-175,240){\line(0,-1){70}}
\put(-105,170){\line(0,-1){70}}
\put(-105,240){\line(0,-1){70}}
\put(-35,170){\line(0,-1){70}}
\put(-35,240){\line(0,-1){70}}
\put(175,170){\line(0,-1){70}}
\put(175,240){\line(0,-1){70}}
\put(105,170){\line(0,-1){70}}
\put(105,240){\line(0,-1){70}}
\put(35,170){\line(0,-1){70}}
\put(35,240){\line(0,-1){70}}
\put(-175,100){\line(0,-1){15}}
\put(-105,100){\line(0,-1){15}}
\put(-35,100){\line(0,-1){15}}
\put(175,100){\line(0,-1){15}}
\put(105,100){\line(0,-1){15}}
\put(35,100){\line(0,-1){15}}
\put(-175,240){\line(0,1){15}}
\put(-105,240){\line(0,1){15}}
\put(-35,240){\line(0,1){15}}
\put(175,240){\line(0,1){15}}
\put(105,240){\line(0,1){15}}
\put(35,240){\line(0,1){15}}
\put(-175,170){\line(1,-1){70}}
\put(-175,170){\line(2,-1){140}}
\put(-175,170){\line(3,-1){210}}
\put(-175,170){\line(4,-1){280}}
\put(-175,240){\line(1,-1){70}}
\put(-175,240){\line(2,-1){140}}
\put(-175,240){\line(3,-1){210}}
\put(-175,240){\line(4,-1){280}}
\put(175,170){\line(-1,-1){70}}
\put(175,170){\line(-2,-1){140}}
\put(175,170){\line(-3,-1){210}}
\put(175,170){\line(-4,-1){280}}
\put(175,240){\line(-1,-1){70}}
\put(175,240){\line(-2,-1){140}}
\put(175,240){\line(-3,-1){210}}
\put(175,240){\line(-4,-1){280}}
\put(-105,170){\line(1,-1){70}}
\put(-105,170){\line(2,-1){140}}
\put(-105,240){\line(1,-1){70}}
\put(-105,240){\line(2,-1){140}}
\put(105,170){\line(-1,-1){70}}
\put(105,170){\line(-2,-1){140}}
\put(105,240){\line(-1,-1){70}}
\put(105,240){\line(-2,-1){140}}
\put(-175,100){\line(1,-1){15}}
\put(-175,100){\line(2,-1){30}}
\put(-175,100){\line(3,-1){45}}
\put(-175,100){\line(4,-1){60}}
\put(175,100){\line(-1,-1){15}}
\put(175,100){\line(-2,-1){30}}
\put(175,100){\line(-3,-1){45}}
\put(175,100){\line(-4,-1){60}}
\put(-105,100){\line(1,-1){15}}
\put(-105,100){\line(2,-1){30}}
\put(105,100){\line(-1,-1){15}}
\put(105,100){\line(-2,-1){30}}
\put(-35,240){\line(-1,1){15}}
\put(-35,240){\line(-2,1){30}}
\put(-35,240){\line(2,1){30}}
\put(-35,240){\line(3,1){45}}
\put(35,240){\line(1,1){15}}
\put(35,240){\line(2,1){30}}
\put(35,240){\line(-3,1){45}}
\put(35,240){\line(-2,1){30}}
\put(-105,240){\line(-1,1){15}}
\put(-105,240){\line(4,1){60}}
\put(105,240){\line(1,1){15}}
\put(105,240){\line(-4,1){60}}
\put(0,50){$\vdots$}
\put(0,290){$\vdots$}
\end{picture}
\end{center}
\caption{\label{fi:sphbra}} 
\centerline {\footnotesize{The stable part of the Bratteli diagram
for the sphere poset $P_6(S^2)$.}}
\vskip1cm
\end{figure}
The set $\{0\}$ is not an ideal. The inductive limit is a subalgebra of bounded
operators on the Hilbert space $\ch = \ch_1 \oplus \cdots \oplus
\ch_8$ with $\ch_i 
,~ i = 1, \dots, 8$ infinite dimensional Hilbert spaces, given by:
\bea
\ca &=& \IC \unit_{\ch_1 \otimes (\ch_5 \oplus \ch_6) \oplus
                          \ch_2 \otimes (\ch_7 \oplus \ch_8)} 
\oplus \IC \unit_{\ch_3 \otimes (\ch_5 \oplus \ch_6) \oplus
                          \ch_4 \otimes (\ch_7 \oplus \ch_8)}
\nonumber \\
& \oplus&
\left( \ck_{\ch_1 \oplus \ch_3} \otimes \IC\unit_{\ch_5 \oplus \ch_6} \right)
\oplus
\left( \ck_{\ch_2 \oplus \ch_4} \otimes \IC\unit_{\ch_7 \oplus \ch_8} \right)
\nonumber \\
& \oplus &
\ck_{\ch_5 \otimes (\ch_1 \oplus \ch_3) \oplus
                         \ch_7 \otimes (\ch_2 \oplus \ch_4)}
\oplus
\ck_{\ch_6 \otimes (\ch_1 \oplus \ch_3) \oplus
                         \ch_8 \otimes (\ch_2 \oplus \ch_4)}
 \nonumber \\
 &~ &\label{sphalg}
\eea

\section{The Behncke-Leptin Construction}\label{se:bl}

Given a poset $P$, there is always an AF algebra $\ca$ such that
$\hat{\ca} = P$. 
A particular procedure to find such an algebra has been described in
the previous 
section, but it is known that there exists more than one \cstar $\ca$ whose
structure space is $P$. For example, if $P$ consists of a single point, we can
take for $\ca$ any of the \cstars $\IM(n,\IC)$ of all $n\times n$ matrices
valued in $\IC$.

It is natural then to ask what are all the algebras associated to a given
finite $T_0$ topological space $P$. This problem was
solved by Behncke and Leptin in 1973. In \cite{BL} they give a complete
classification of all separable \cstars $\ca$ with finite structure
spaces. Such classification requires the definition of a function $d$ on $P$,
called defector, valued in $\bar \IN=\{\infty, 0,1,2,...\}$. Given $P$ and
$d$, the Behncke-Leptin construction gives 
a separable \cstar $\ca(P,d)$ such that $\hat \ca(P,d)=P$. Furthermore, any
separable \cstar $\ca$ satisfying $\hat \ca =P $  is isomorphic to
$\ca(P,d)$ for some $d$.

A {\it defector }
$d$ on the poset $P$ is a $\bar \IN$-valued function on $P$ such that
\be
d(x)>0 ~~\mbox{ if $x$ is maximal}. \label{5.1}
\ee
Two defectors $d$ and $d'$ are declared to be equal if there exists an
automorphism $\varphi $ of $P$ such that $d'=d\circ \varphi $.
They are called immediately equivalent if
$d(x)=d(x')$ for all $x\in P$ with the exception of at most one
nonmaximal $y\in P$, such that
\[
d(y)=d'(y)+d'(z)~~~{\rm or}~~~d'(y)=d(y)+d(z)
\]
for some $z$ covering $y$ \footnote{We say that $y$ covers $x$ if
$x\prec y$ and 
there is no $z$ such that $x \prec z \prec y$.}, if $d(z)=d'(z)<\infty
$. In the 
case $d(z)=\infty$, $d(y)$ and $d'(y)$ may be arbitrary. Then two defectors
are defined to be equivalent ($d\sim d'$) if there exists a finite
sequence of immediately equivalent defector connecting them.

We will start by describing the Behncke-Leptin construction for a
special class of posets called forests. Then we will give the
generalization for 
an arbitrary finite poset.

\subsection{The Behncke-Leptin Construction for a Forest}\label{sub:5.1}

A forest is a poset $F$ such that
\be
\{x,y,z\in F, x\preceq z, y\preceq z\}\Rightarrow \{x\preceq y
\mbox{ or } y\preceq x \}. \label{5.2}
\ee

Given a forest $F$ and a defector $d$ on $F$, the Behncke-Leptin construction
consists of the following steps. First we introduce a Hilbert
space $\ch(F,d)$ associated to the whole forest $F$. Second, for each point
$x\in F$, we introduce a subspace $\ch(x)\subseteq \ch(F,d)$ and a set
of operators 
$\czr_x$ acting on $\ch(x)$. Actually, $\czr_x$ can be thought of as acting on
the whole $\ch(F,d)$ by defining its action on the complement of $\ch(x)$ to be
zero. Then the \cstar $\ca$ associated to the forest is the one
generated by the 
$\czr_x$'s as $x$ varies in $F$.

Now we explain how to determine $\ch(F,d), \ch(x)$ and $\czr_x$. The Hilbert
space $\ch(F,d)$ can be obtained using an auxiliary forest $F'$
constructed from 
$F$ in the following way. The forest $F'$ contains a point
$x^{(1)}$ for each maximal point $x\in F$ and a pair of points $x^{(1)}_i$ and
$x^{(2)}_i$  for each non maximal point $x_i\in F$. Then on $F'$ we introduce a
partial order by declaring that $x_i^{(2)}$ is covered by both $x_j^{(1)}$ and
$x_j^{(2)}$ if and only if $x_i$ is covered by $x_j$. Figure
\ref{fi:forest} shows an example of $F$ and the corresponding $F'$.

\begin{figure}
\begin{center}

\begin{picture}(400,150)(-20,-30)

\put(100,50){\circle*{4}}
\put(50,100){\circle*{4}}
\put(150,100){\circle*{4}}
\put(100,0){\circle*{4}}

\put(100,50){\line(-1,1){50}}
\put(100,50){\line(1,1){50}}
\put(100,50){\line(0,-1){50}}

\put(105,47){$x_2$}
\put(45,105){$x_3$}
\put(145,105){$x_4$}
\put(105,-2){$x_1$}
\put(96,-18){$F$}

\put(250,50){\circle*{4}}
\put(200,100){\circle*{4}}
\put(300,100){\circle*{4}}
\put(250,0){\circle*{4}}
\put(225,50){\circle*{4}}
\put(225,0){\circle*{4}}

\put(250,50){\line(-1,1){50}}
\put(250,50){\line(1,1){50}}
\put(250,50){\line(0,-1){50}}
\put(250,0){\line(-1,2){25}}

\put(255,47){$x_2^{(2)}$}
\put(195,105){$x_3^{(1)}$}
\put(295,105){$x_4^{(1)}$}
\put(255,-2){$x_1^{(2)}$}
\put(205,-2){$x_1^{(1)}$}
\put(205,47){$x_2^{(1)}$}
\put(246,-18){$F'$}

\end{picture}

\end{center}
\caption{\label{fi:forest}}
\centerline{{\footnotesize An example of a forest $F$ and the auxiliary
forest $F'$.}}
\vskip1cm
\end{figure}

In $F'$ we consider the maximal chains
$C_\a=\{x^{(p_1)}_ 1,x^{(p_2)}_2,....,x^{(p_k)}_k\}$, which can be seen to be
necessarily of the form
\be
C_\a=\{x^{(2)}_1,x^{(2)}_2,....,x^{(2)}_{k-1},x^{(1)}_k\}\; . \label{5.3}
\ee
For example, in $F'$ of Figure \ref{fi:forest}, the maximal chains are
$\{x_1^{(1)}\},\{x_1^{(2)},x_2^{(1)}\},\{x_1^{(2)},x_2^{(2)},x_3^{(1)}\},
\{x_1^{(2)},x_2^{(2)},x_4^{(1)}\}$.

To each maximal chain $C_\a$ we associate the Hilbert space
\be
h(C_\a)=l_{x_1}\otimes l_{x_2}\otimes ...\otimes l_{x_{k-1}}\otimes {\bf
C}^{d(x_k)}\; , \label{5.4}
\ee
where $d(x_k)$ is the value of the defector $d$ at the point
$x_k\in F$ and $l_{x_i}$ can be realized as the Hilbert space $\ell^2$
of sequences 
$(f_1,f_2,....)$ of complex numbers with $\sum_n |f_n|^2 < \infty$.
We then define the total Hilbert space $\ch(F,d)$ associated to $F$ to be
\be
\ch(F,d)= \bigoplus_\a h(C_\a) \label{5.5} \; ,
\ee
where we sum over all maximal chains $C_\a $ in $F'$.

In a similar
way, we introduce the subspaces $\ch(x_i)$ associated to a single point $x_i
\in F$ by
\be
\ch(x_i)=\bigoplus _\b h(C_\b) ~~\mbox{for all $C_\b$ such that
$x_i^{(p)}\in C_\b $}\; .\label{5.6}
\ee
Notice that if we consider the subforest $F_x$ of $F$ given by
\be
F_x=\{y\in F|x\preceq y\}\; ,\label{5.7}
\ee
we can construct the Hilbert space $\ch(F_x,d_x)$,
where $d_x$ is the restriction of $d$ to $F_x$. An important property
of $\ch(x)$ defined in (\ref{5.6}) is that it satisfies
\be
\ch(x)=\ch_{x}\otimes \ch(F_x,d_x) \label{5.8}
\ee
where
\be
\ch_{x}=\bigotimes _i l_{x_i}\, ~~ \mbox{for all $x_i\prec x$} \label{5.9}
\ee
and $\ch_x = \IC$ if $x$ is a minimal point.

Now we are ready to define the \cstar $\ca(F,d)$. First, let us
introduce the algebra of operators $\czr_x$, acting on $\ch(x)$, given by
\be
\czr_x=\IC \unit_{\ch_x}\otimes \ck(\ch(F_x,d_x))\label{5.10}\; ,
\ee
$\unit_{\ch_x}$ being the identity operator on $\ch_x$ and
$\ck(\ch(F_x,d_x))$ being the algebra of compact operators on $\ch(F_x,d_x)$.
In other words, $\czr_x$ acts as multiples of the identity on the Hilbert space
$\ch_x$ determined by the points $x_i \prec x$  which precede $x$, as
in (\ref{5.9}), and as compact operators on the Hilbert space $\ch(F_x,d_x)$
determined by the points $x_j\succeq x$ which follow $x$.

Then
$\ca(F,d)$ is the algebra of operators on $\ch(F,d)$ generated by all
$\czr_x$ as 
$x$ varies in $F$.

The algebras $\czr_x$, with $x\in F$, satisfy:
\[
\czr_x\czr_y\subset \czr_x ~~~\mbox{if $x\preceq y$}
\]
and
\be
\czr_x\czr_y=0 ~~~\mbox{if $x$ and $y$ are incomparable}. \label{5.11}
\ee

One of the major  results of \cite{BL} is the following
theorem, which establishes that the structure space of the \cstar
$\ca(F,d)$ constructed according to the rules given above is homeomorphic to
the forest $F$:
\bprop
Let $F$ be a finite forest with defector $d$ and $\ca(F,d)$ the
algebra of operators on $\ch(F,d)$ defined as above. Then we have:
\begin{itemize}
\item[(i)] if $E$ is a closed subset of $F$ with complement $U$, then
$I_E=\bigotimes _{x\in U}\czr_x$ is a closed two-sided ideal of
$\ca(F,d)$, and $A_E=\bigotimes _{x\in E}\czr_x$ is a closed subalgebra of
$\ca(F,d)$;
\item[(ii)] every two-sided ideal of $\ca(F,d)$ is of the form $I_E$
for some closed 
$E\subset F$ and $I_E$ is primitive iff $E=\overline{\{x\}}$. In
particular, $\hat \ca(F,d)=F$.
\end{itemize}
\eprop

Let us illustrate the Behncke-Leptin construction for a very simple
forest, namely the $\bigvee$ poset of Figure \ref{fi:vposet}. The correspondent
associated forest $P'$ is illustrated in Figure  \ref{fi:vforest}.

\begin{figure}[htb]
\begin{center}

\begin{picture}(320,150)(-60,20)

\put(100,50){\circle*{4}}
\put(50,100){\circle*{4}}
\put(150,100){\circle*{4}}
\put(75,50){\circle*{4}}

\put(100,50){\line(-1,1){50}}
\put(100,50){\line(1,1){50}}

\put(98,35){$q^{(2)}$}
\put(45,105){$p_1^{(1)}$}
\put(145,105){$p_2^{(1)}$}
\put(70,35){$q^{(1)}$}

\end{picture}

\end{center}
\caption{\label{fi:vforest}}
\centerline{{\footnotesize  The forest associated to the $\bigvee$ poset.}}
\vskip1cm
\end{figure}

We consider a generic defector $d$. From the diagram of
$P'$ in Figure \ref{fi:vforest} we can write down all its maximal chains:
\[ \{q^{(2)},p_1^{(1)}\} \;, \; \{q^{(2)},p_2^{(1)}\} \; , \;
\{q^{(1)}\} \]
and following (\ref{5.4}) and (\ref{5.5}) we see that
$\ch(F,d)$ is given by
\be
\ch(F,d)=\left( l_{q} \otimes \IC^{d(p_1)} \right) \oplus
       \left( l_{q} \otimes \IC^{d(p_2)} \right) \oplus \IC^{d(q)}
\; .\label{5.12}
\ee
The subspaces $\ch(x_i)$ can also be determined from the diagram
of $P'$:
\begin{eqnarray}
 \ch(p_1) & = & l_{q} \otimes \IC^{d(p_1)} \nonumber \\
 \ch(p_2) & = & l_{q} \otimes \IC^{d(p_2)} \nonumber \\
 \ch(q) & = & \ch(F,d) \; . \label{5.13}
\end{eqnarray}
Notice that the factorization expressed in (\ref{5.8}) is satisfied,
where now
\begin{eqnarray}
\ch_{p_1}  = l_{q}  & , & \ch(F_{p_1},d_{p_1}) = \IC^{d(p_1)} \nonumber \\
\ch_{p_2}  = l_{q}  & , & \ch(F_{p_2},d_{p_2}) = \IC^{d(p_2)} \nonumber \\
\ch_{q}    = \IC& , & \ch(F_q,d_q) = \ch(F,d) \, . \label{5.14}
\end{eqnarray}

The \cstar $\ca(F,d)$ is generated by all $\czr_x$, $x\in F$. The latter
reads
\begin{eqnarray}
\czr_{p_1} & = & \IC \unit_{\ch_{p_1}}\otimes \ck(\IC^{d(p_1)}) \nonumber \\
\czr_{p_2} & = & \IC \unit_{\ch_{p_2}}\otimes \ck(\IC^{d(p_2)}) \nonumber \\
\czr_{q} & = & \ck(\ch(F,d))\; . \label{5.15}
\end{eqnarray}

Notice that for the defector $d(p_1)=d(p_2)=1$ and $d(q)=0$ we get $\ch(F,d) =
\ch_1 \oplus \ch_2$ and $\ca = \IC \unit_1 + \IC \unit_2 + \ck_{12}$
and thus recover the algebra we got for the $\bigvee$ poset via the Bratteli
construction in section \ref{sub:dia}.

\subsection{The Behncke-Leptin Construction for Posets}\label{sub:5.2}

To generalize the procedure of the last section to an arbitrary
poset $P$ with defector $d$, we have first to introduce  a
forest $\bar P$, uniquely determined by $P$.

Let $P$ be a finite poset. A rope $r$ of $P$ is a (not
necessarily maximal) chain in $P$ starting
from a minimal element and ending at some $x\in P$. The set $\bar P$ of all
ropes of $P$ ordered by inclusion is a poset. One can show that $\bar
P$ is in fact a forest. Let $\varphi :\bar P\rightarrow P$ denote the
surjective map which assigns to each rope $r\in \bar P$ its end point
$\varphi (r)\in P$. Following \cite{BL}, we will call the pair $(\bar
P,\varphi )$ the covering forest of $P$. An example
is given in Figure \ref{fi:circleforest}, which shows the covering
forest of the 
circle poset $P_4(S^1)$ of Figure \ref{fi:hassecircle}.

\begin{figure}
\begin{center}

\begin{picture}(400,150)(0,20)

\put(100,50){\circle*{4}}
\put(50,100){\circle*{4}}
\put(150,100){\circle*{4}}

\put(100,50){\line(-1,1){50}}
\put(100,50){\line(1,1){50}}

\put(91,37){$\{x_1\}$}
\put(33,105){$\{x_1,x_2\}$}
\put(133,105){$\{x_1,x_4\}$}

\put(300,100){\circle*{4}}
\put(200,100){\circle*{4}}
\put(250,50){\circle*{4}}

\put(250,50){\line(-1,1){50}}
\put(250,50){\line(1,1){50}}

\put(241,37){$\{x_3\}$}
\put(183,105){$\{x_3,x_2\}$}
\put(283,105){$\{x_3,x_4\}$}

\end{picture}

\end{center}
\caption{\label{fi:circleforest}}
\centerline{{\footnotesize The covering forest of the circle poset
$P_4(S^1)$.}}
\vskip1cm
\end{figure}

Given a defector $d$ on $P$ we define a defector $\bar d$ on $\bar P$
in a natural way via the pull-back:
\be
\bar d=d\circ \varphi \; .\label{5.16}
\ee

Then, since $\bar P$ is a forest, we can construct the algebra $\ca(\bar
P,\bar d)$ following section \ref{sub:5.1}.

Finally, to identify the \cstar $\ca(P,d)$ associated to the poset $P$ and the
defector $d$, we proceed to realize $\ca(P,d)$ as a subalgebra of $\ca(\bar
P,\bar d)$.  In order to do so, we need to point out a simple property of the
covering forest $(\bar P,\varphi )$.

Let $r$, $s\in \bar P$ be in the inverse image $\varphi^{-1}(x)$ of
$x\in P$. Then, the subforest $(\bar P)_r$ (see (\ref{5.7})) is
naturally isomorphic to $(\bar P)_s$. Indeed, $(\bar P)_r$ and $(\bar
P)_s$ consist of all extensions of the rope $r$ and $s$ respectively. By
hypothesis, $r$ and $s$ have the same end point $x\in P$, so that
$(\bar P)_r\sim (\bar P)_s$. Thus
\[\ck\left( \ch((\bar P)_r,\bar d_r) \right) \simeq
\ck\left( \ch((\bar P)_s,\bar d_s) \right) \equiv \ck_x \; ,\]
so that the algebras $\czr_s$, $\czr_r\in \ca(\bar P,\bar d)$ are given by
\[ \czr_r=\IC \unit_{\ch_r}\otimes \ck_x \; , \]
\[ \czr_s=\IC \unit_{\ch_s}\otimes \ck_x \; . \]

For each $x\in P$ we define the algebra $A_x$
\be
A_x= \bigoplus _{r\in \varphi ^{-1}(x)}\czr_r
\ee
and a subalgebra $\czr_x\subset A_x$ given by all elements $a\in A_x$
of the form 
\[
  a=(\l_{r_1} \unit_{\ch_{r_1}}\otimes k)+(\l_{r_2}
\unit_{\ch_{r_2}}\otimes k)+ ...  +
(\l_{r_n} \unit_{\ch_{r_n}}\otimes k)\; ,
\]
where $r_i \in \varphi ^{-1}(x)$, $\l_{r_j} \in \IC$ and $k \in \ck_x$.
Thus
\be
\czr_x= \{a\in A_x\,|\, a=\bigoplus _{r\in \varphi ^{-1}(x)} (\l_{r}
\unit_{\ch_r}\otimes k) \mbox{ , } \l_r \in \IC \mbox{ and } k \in \ck_x \} \;
.\label{5.18}
\ee
The \cstar $\ca(P,d)$ that satisfies
\be
\hat \ca(P,d)=P
\ee
is then generated by all $\czr_x$ with $x\in P$.

There is an intuitive interpretation for (\ref{5.18}). The poset $P$
can be obtained from $\bar P$ by identifying any two ropes $r$ and
$s$ that have the same ending point. Equation (\ref{5.18}) simply expresses
this identification at an algebraic level.

For example, for the circle poset $P_4(S^1)$ these rules give the following
algebras, acting on $\ch = \ch_1 \oplus \ch_2 \oplus \ch_3 \oplus \ch_4$:
\bea
\ca_{x_4} &=& \IC \unit_1 + \IC \unit_3 \nonumber \\
\ca_{x_2} &=& \IC \unit_2 + \IC \unit_4 \nonumber \\
\ca_{x_1} &=& \IC \unit_1 + \IC \unit_2 + \ck_{12} \nonumber \\
\ca_{x_3} &=& \IC \unit_3 + \IC \unit_4 + \ck_{34}  \; ,
\eea
if one chooses the defector $d(x_1)=d(x_2)=1$, $d(x_3)=d(x_4)=0$. Thus the
algebra associated to $P_4(S^1)$ is
\be
\ca = \IC \unit_1 + \IC \unit_2 + \IC \unit_3 + \IC \unit_4 +
\ck_{12} + \ck_{34} \; . \label{algex}
\ee
As before this is the algebra one gets  for $P_4(S^1)$ by means of the
Bratteli construction explained in section \ref{sub:pobra}. 

Equivalent defectors give rise to isomorphic $C^*$-algebras, whereas
by choosing 
different non-equivalent defectors one can construct non-isomorphic
\cstars that 
all have $P$ as structure space. In this way one can obtain {\it all} \cstars
$\ca$ whose (finite) dual $\hat{\ca}$ is homeomorphic to the poset
$P$, as it is 
established in \cite{BL}, which
we quote in conclusion of this section:

\bprop
\begin{itemize}
\item[(i)] Every separable \cstar $\ca$ with finite dual $\hat \ca =P$
is isomorphic to some $\ca(P,d)$.
\item[(ii)] $\ca(P,d)$ is isomorphic to $\ca(P,d')$ if and only if $d$
and $d'$ are equivalent.
\end{itemize}
\eprop

\section{Final remarks}

In this article, we have seen how a finite poset is truly a
``noncommutative space" 
or ``noncommutative lattice", since it can be described as the
structure space of 
a noncommutative \cstar $\ca$, which turns out to be always a postliminal AF
algebra. We have also seen as this correspondence is not one-to-one,
more than one 
non-isomorphic \cstar leading to the same poset. This relation between
posets and 
\cstars was used in \cite{pangs} to give a dualization of the
approximation method 
for topological spaces introduced in \cite{So}.

In our previous work \cite{lisbon} we have showed how it is possible
to construct a 
quantum theory on posets, by making use of the corresponding $C^*$-algebra. We
have also seen how important topological properties of the continuum, such as
homotopy, can be captured by the poset approximation and manifest themselves in
the corresponding quantum mechanics.

We are thus naturally led to examine how one can construct further geometric
structures on posets, as it suggested by Connes' noncommutative geometry
\cite{Co}. First of all, we are interested in the construction of bundles and
characteristic classes over a poset and, as a first step in this direction, one
should examine the K-theory of these noncommutative lattices. This is the topic
discussed in  \cite{kt}, where we present a study of the algebraic
K-theory of AF 
algebras associated to a poset. 

Then one would like to construct bundles, and notably nontrivial
bundles, over a 
poset, and consider, for instance, the analogue of the monopole bundle over the
lattice approximating the two-dimensional sphere and of nontrivial
``topological 
charges". Work in this direction is in progress.

\vskip1cm
{\Large {\bf Acknowledgments}}

This work was initiated while the authors were at Syracuse University. 
We thank A.P. Balachandran, G. Bimonte, F. Lizzi e G. Sparano for many
fruitful discussions and useful advice. The final version was written while
G.L. and P.T-S were at ESI in Vienna. They would like to thank G. Marmo and P.
Michor for the invitation and all people at the Institute for the warm
and friendly 
atmosphere.

We thank the `Istituto Italiano per gli Studi Filosofici' in Napoli for
partial support.
The work of P.T-S. was also supported by the Department of Energy,
U.S.A. under
contract number DE-FG-02-84ER40173. The work of G.L. was
partially supported by the Italian `Ministero dell' Universit\`a e
della Ricerca Scientifica'.

\end{document}